\documentclass{amsart}
\usepackage[margin=1.25in,dvips]{geometry}

\usepackage{cite}
\usepackage{amssymb,mathrsfs,graphicx}
\usepackage{amsmath}
\usepackage{bbm}
\usepackage{xcolor}
\usepackage{soul}

\numberwithin{equation}{section}

\title[Local Decay Estimates]{Local Decay Estimates}

 \author[Avy Soffer]{Avy Soffer}
 \address[Avy Soffer]{\newline
        Department of Mathematics, \newline
         Rutgers University, New Brunswick, NJ 08903 USA.}
  \email[]{soffer@math.rutgers.edu}

\author[Xiaoxu Wu]{Xiaoxu Wu}
\address[Xiaoxu Wu]{\newline
        Mathematical Sicences Institute, \newline
        Australia National University, Acton, ACT 2601, Australia}
 \email[]{Xiaoxu.Wu@anu.edu.au}

\newtheorem{theorem}{Theorem}[section]
\newtheorem{lemma}{Lemma}[section]

\newtheorem{assumption}{Assumption}[section]
\newtheorem{corollary}{Corollary}[section]
\newtheorem{proposition}{Proposition}[section]
\newtheorem{remark}{Remark}[section]

\newtheorem{definition}{Definition}[section]

\newcommand{\s}{ \mathcal{L}}
\usepackage{xcolor}

\newcommand{\eq}{\begin{equation}}
\newcommand{\eeq}{\end{equation}}

\begin{document}

\date{\today}



\begin{abstract}
\noindent We give a proof of local decay estimates for Schr\"odinger-type equations, which is based on the knowledge of Asymptotic Completeness (AC).
This approach extends to time dependent potential perturbations, as it does not rely on Resolvent Estimates or related methods. Global in time Strichartz estimates follow for quasi-periodic time-dependent potentials from our results.

\end{abstract}
\maketitle
\centerline{\date}

\tableofcontents

\section{Introduction}\label{sec:intro}

\noindent Local decay estimates are a priori estimates of the solutions of dispersive equations.
It states that the solution (for an initial condition associated with scattering) decays at least in an integrable rate in time, in every compact region of space, for a dense set, in the scattering subspace (see Eq.~\eqref{def: scattering}), of initial conditions. An equivalent statement is that the resolvent of the Hamiltonian of the dynamics is bounded on properly weighted $L^2$ space.

Such estimates play a crucial role in scattering theory, as they imply the existence and completeness of the M\"{o}ller wave operators. Moreover, local decay and other propagation estimates have important applications beyond the proof of asymptotic completeness (AC). For example, they are essential in linear and nonlinear time-dependent resonance theory \cite{S-Wei2, S-Wei1}. Notable consequences of local decay include Strichartz estimates and the propagation of regularity for nonlinear dispersive equations \cite{T2006}. These estimates provide insights into the rate of convergence with asymptotic dynamics.

Proving such estimates for an interacting system is indeed challenging, as we do not have direct application of the method of stationary phase. However, it is possible to prove asymptotic completeness (AC) without using local decay. This was first shown by Enss in 1978 \cite{E1978}. The Enss method was applicable to both two-particle and three-particle scattering in quantum mechanics. However, proving AC for four or more particles required the use of local decay \cite{SS1987}. This was achieved using the Mourre estimate combined with Mourre's method of differential inequalities \cite{M1979, M1981, PSS, FH}, or the method of propagation estimates developed by Sigal and Soffer \cite{SS1988}. See also \cite{Skib, HSS1999, Ger, G1990}.

 Recently, Liu and the first author introduced a general approach to proving asymptotic completeness (AC) \cite{Liu-Sof1, Liu-Sof2}. Their method applies to both linear and nonlinear dispersive equations and requires localization assumptions on the interaction terms. Specifically, their approach necessitates radial symmetry for nonlinear equations.

Subsequently, the authors improved Liu and Soffer's method by constructing the free channel wave operator in an adapted way \cite{SW20221, Sof-W2, Sof-W3}. This enhanced method does not require radial symmetry or localization assumptions on the interaction terms. It is applicable to both linear and nonlinear dispersive equations, including those with time-dependent interaction terms and Klein-Gordon type equations.

  Neither method directly employs local decay. However, proving local decay in these cases is not straightforward. In this work, we use the knowledge of asymptotic completeness (AC) obtained from \cite{SW20221} to derive local decay estimates. We first establish this for the Schrödinger equation with a time-independent potential that is localized in space, and then for a potential that is quasi-periodic in time and also localized in space.

Time-dependent potentials have been studied previously. In particular, Rodnianski and Schlag \cite{RS2004} established \(L^p\) estimates for time-dependent potentials under the assumption of smallness on the size of the potential. Recent progress in this direction can also be found in \cite{G2009, V2013, SW2020, GH, B2011,EG2010}.

\subsection{Main results}
\subsubsection{The time-dependent result} Let $H_0:=-\Delta_x$ and $\langle \cdot\rangle: \mathbb{R}^n\to \mathbb{R}, x\mapsto \sqrt{|x|^2+1}, n\geq 1$. We consider the Schr\"odinger equation with a time-dependent potential $V(x,t)$ 
\eq
\begin{cases}
i\partial_t\psi(x,t)=(H_0+V(x,t))\psi(x,t)\\
\psi(x,0)=\psi_0\in \s^2_x(\mathbb{R}^n)
\end{cases},\qquad (x,t)\in \mathbb{R}^n\times \mathbb{R},\quad n\geq 5,\label{tSE}
\eeq
where $V(x,t)$ satisfies 
\begin{assumption}\label{aspV2} $V(x,t)$ is real-valued and satisfies the condition $\langle x\rangle^\sigma V(x,t) \in L^\infty_{x,t}(\mathbb{R}^{n+1})$ for some $\sigma > 6$, where $n \geq 5$. Additionally, $V(x,t)$ can be expressed as
\begin{equation}
V(x,t) = \sum_{j=1}^N V_j(x,t), \qquad N \in \mathbb{N}^+,
\end{equation}
where each $V_j(x,t)$, for $j = 1, \dots, N$, is continuous and periodic in $t$ in the sense of 
\begin{equation}
    \lim\limits_{s\to t} \| V_j(x,t)-V_j(x,s)\|_{\mathcal L^\infty_x(\mathbb R^n)}=0.\label{asp: conti}
\end{equation}
\end{assumption}
\noindent This potential $V$ is quasiperiodic in $t$. The definition of a quasi-periodic function is provided below.
\begin{definition}[Quasi-periodic functions]\label{def: quasif} A function $f(t)$ is \emph{quasi-periodic} if there exist $N$ tori $\mathbb T_j:=[0, T_j), j=1,\cdots,N (N\in \mathbb N^+)$ and a function $g: \mathbb T_1\times \cdots \times \mathbb T_N\to \mathbb C$ such that whenever 
\begin{equation}
    s_j=t\, \quad\text{mod}\,\, T_j, \qquad j=1,\cdots,N,
\end{equation}
\begin{equation}
    f(t)=g(s_1,\cdots,s_N).
\end{equation}
We say that two quasi-periodic functions are of the same type if the domains of their corresponding functions $g$ are the same. 

\end{definition}
\noindent An example of quasi-periodic functions which are not periodic is $f(t)=\sin(t)+\sin(\sqrt{2}t)$, where the periods of $\sin(t)$ and $\sin(\sqrt{2}t)$ are $2\pi$ and $\sqrt{2}\pi$ respectively. Since the ratio of these periods is $\frac{2\pi}{\sqrt{2}\pi}=\sqrt{2}$, which is irrational, the function does not repeat itself over any common period, making it quasi-periodic.

Throughout this paper, $H^k(\mathbb R^n),$ for $ k\geq 0$, denotes the Sobolev space and $A\lesssim_a B$ means $A\leq CB $ for some constant $C=C(a)>0$. We also use the notation 
\begin{equation}
    \|\cdot\|
\end{equation}
to denote either
\begin{equation}
\|\cdot\|_{\mathcal L^2_x(\mathbb{R}^n)}\quad\text{or}\quad\|\cdot\|_{\mathcal L^2_x(\mathbb{R}^n)\to \mathcal L^2_x(\mathbb R^n)}
\end{equation}
depending on the context. Sometimes we use $\|\cdot\|_{2\to 2}$ to denote $\|\cdot\|_{\mathcal L^2_x(\mathbb R^n)\to \mathcal L^2_x(\mathbb R^n)}$.

Under Assumption~\ref{aspV2}, $H\equiv H_0+V(x,t)$ is self-adjoint on the domain of $H_0$, i.e. $V(x,t)$ is $H_0$-bounded with the relative bound $<1$: for all $u\in H^2_x(\mathbb R^n)$,
\begin{equation}
    \|Vu\|\leq a\| H_0 u\|+b\|u\|, \qquad\exists \,0\leq a<1,\, b>0,\label{sa. ineq}
\end{equation}
and the evolution operators $U(t,s), t,s\in \mathbb R,$ generated by $H$ is unitary on $\mathcal L^2_x(\mathbb R^n)$. Indeed, for all initial data $\psi_0\in \mathcal L^2_x(\mathbb R^n)$, the solution $\psi(x,t)$ exists in $\mathcal L_x^2(\mathbb R^n)$ for all $t\in \mathbb R$. See e.g. Theorem X.70 of \cite{RSII}.

In this note, we stick to the case when $t\geq 0$. We refer to the initial state $\psi_0\in \mathcal L^2_x(\mathbb R^n)$ as a \emph{scattering state} if $\psi(x,t)$ evolves like a free flow asymptotically, that is, for some $\psi_+\in \mathcal L^2_x(\mathbb R^n)$,
\begin{equation}\label{def: scattering}
    \|\psi(x,t)-e^{-itH_0}\psi_+\|\to 0,\qquad \text{ as }t\to \infty.
\end{equation}
The space of all scattering states is a subspace of $\mathcal L^2_x(\mathbb R^n)$ and such space varies as the initial time $t_0$ varies. We introduce the projection on the space of all scattering states at the initial time $t_0$ below. For $m>0$, let $F(k\geq m)\equiv F(\frac{k}{m}) (m>0)$ denote a smooth characteristic function with
\eq
F(\lambda)=\begin{cases}1&\text{ when }\lambda>1\\
0& \text{ when }\lambda<1/2\end{cases},
\eeq
and let $F(k<m)=\mathbbm1-F(k\geq m)$ denote the complement of $F(k\geq m)$. The projection on the space of all scattering states is given by, with $p:=-i\nabla_x$, for $\alpha\in (0, 1-\frac{2}{n}), n\geq 3$ and $t_0\in \mathbb R$,
\begin{equation}
    P_c(t_0):= s\text{-}\lim\limits_{t\to \infty} U(t_0,t_0+t) F(\frac{|x-2tp|}{t^\alpha}< 1) U(t_0+t,t_0),\qquad \text{ on }\mathcal L_x^2(\mathbb R^n), \label{def: Pct0}
\end{equation}
which is constructed in~\cite{SW20221} provided that $V(x,t)\in \mathcal L^\infty_t\mathcal L^2_x(\mathbb R^{n+1}), n\geq 3$. It is worth noting that the r.h.s. of Eq.~\eqref{def: Pct0} is independent of $\alpha$ (See Eq.~\eqref{WO: def} below and Eq.~(2.13) of \cite{SW20221}) and $P_c(t_0)$ defined in Eq.~\eqref{def: Pct0} satisfies $P_c(t_0)=\Omega_+(t)\Omega_+^*(t)$, where 
\begin{equation}
\Omega_+(t):=s\text{-}\lim\limits_{s\to \infty}U(t,t+s)e^{-isH_0} \qquad \text{ on }\mathcal L^2_x(\mathbb R^n) 
\end{equation}
and $\Omega_+^*(t)$ is the conjugate of $\Omega_+(t)$.

Throughout this paper, $C$ will denote a constant and may vary from one line to another. We write $\lesssim$ or $\gtrsim$ whenever $A\leq CB$ or $CA\geq B$ for some constant $C>0$. We write $A\lesssim_a B$ or $A\gtrsim_a B$ if  $A\leq C_aB$ or $C_aA\geq B$ for some constant $C_a>0$ which depends on parameter $a$. The Fourier transform and its inverse are given by 
\begin{equation}
    \hat f(\xi):=\frac{1}{(2\pi)^{n/2}}\int e^{-ix\cdot \xi} f(x)dx
\end{equation}
and 
\begin{equation}
    f(x):=\frac{1}{(2\pi)^{n/2}}\int e^{ix\cdot \xi} \hat f(\xi)d\xi
\end{equation}
for $f\in \mathcal L_x^2(\mathbb R^n)$.

We prove first that, provided Assumption~\ref{aspV2} holds true, there exists $p_0\in [2,\infty)$ such that
\begin{equation}
\left( \int_0^\infty    \| \langle x\rangle^{-\eta} U(t+t_0,t_0)P_c(t_0)\psi_0\|^pdt\right)^{\frac{1}{p}}\lesssim \|\psi_0\|\label{goal: 1}
\end{equation}
for all $\eta>\frac{3}{2}$ and $p_0\leq p<\infty$.
\begin{theorem}\label{thm}If Assumption~\ref{aspV2} is satisfied, then with $p_0=2$, \eqref{goal: 1} holds true for all $p_0\leq p<\infty$, when $n\geq 8$.
\end{theorem}
We refer to estimate~\eqref{goal: 1} as $L^p$ local decay estimate. It is known that the free flow satisfies $L^2$ local decay:
\begin{equation}
  \int_0^\infty \| \langle x\rangle^{-\delta} e^{-itH_0} \psi_0\|^2 dt\lesssim_{1-\delta} \|\psi_0\|^2,\qquad \forall\, \delta>1.\label{free: diseq1}
\end{equation}
So the conclusion of Theorem~\ref{thm} is not sharp. Next, we explain what additional conditions are needed to obtain the $L^2$ local decay estimate
\begin{equation}
 \int_0^\infty    \| \langle x\rangle^{-\eta} U(t+t_0,t_0)P_c(t_0)\psi_0\|^2dt\lesssim\|\psi_0\|^2,\qquad \text{for some }\eta>1.\label{goal: 2}
\end{equation}
We need extra properties of $P_b(t_0)\equiv\mathbbm1-P_c(t_0)$. The existence of $P_c(t_0)$ implies the existence of $P_b(t_0)$. To be precise, we need 
\begin{assumption}\label{aspV3}Either 
\eq
\sup\limits_{t_0\in \mathbb R}\| P_b(t_0)\langle x\rangle^\delta\|<\infty\label{aspeq1}\qquad\text{ for some }\delta>\frac{3}{2}
\eeq
or 
\begin{equation}
    \sup\limits_{t_0\in \mathbb R}\| P_b(t_0)|p|^{-1/2}\langle x\rangle^{1/2+\epsilon}\|<\infty\qquad\text{ for some }\epsilon>0\label{aspeq2}
\end{equation}
holds true.
\end{assumption}
Indeed, in the proof of $L^2$ local decay, we have to estimate $P_b(t)e^{-itH_0}f, f\in \mathcal L^2_x(\mathbb R^n)$. If we write 
\begin{equation*}
    P_b(t)e^{-itH_0}f=\left(P_b(t)\langle x\rangle^\delta \right)\left( \langle x\rangle^{-\delta}e^{-itH_0}f\right)=\left(P_b(t)|p|^{-1/2}\langle x\rangle^{1/2+\epsilon} \right)\left( \langle x\rangle^{-1/2-\epsilon}|p|^{1/2}e^{-itH_0}f\right),
\end{equation*}
then, intuitively, Conditions~\eqref{aspeq1} and~\eqref{aspeq2} are corresponding to $L^2$ local decay estimate and 
\begin{equation}
  ( \text{$L^2$ local smoothing estimate})\qquad   \int_0^\infty \| \langle x\rangle^{-1/2-\epsilon}|p|^{1/2} e^{-itH_0} \psi_0\|^2 dt\lesssim_{\epsilon} \|\psi_0\|^2,\label{free: diseq2}
\end{equation}
respectively. We believe Assumption~\ref{aspV3} is a necessary condition based on our experience with the time-independent potential. When $H=H_0+V$ is time-independent and $V$ satisfies, for example, $\langle x\rangle^{1+\epsilon}|V(x)|\lesssim 1$ for some $\epsilon>0$,
$P_b(t_0)$ is the projection onto the space of the discrete spectrum of $H$. It is known that, additionally, if $0$ is neither an eigenvalue nor a resonance of $H$, $H$ has finitely many eigenvalues, each with finite multiplicity. It is also known that in $5$ or higher space dimensions, there is no resonance. See, for example,~\cite{RS2004} and the references therein. In the case of~\eqref{tSE}, we prove in Proposition~\ref{prop: weight} that for all $\delta\in [0, \frac{n}{2}-2), n\geq 5$,
\begin{equation}
   \sup\limits_{t_0\in \mathbb R} \| P_b(t_0)\langle x\rangle^\delta \|\lesssim_{\delta,n} 1,\label{asp: Pbineq}
\end{equation}
under Assumption~\ref{aspV2}. In the proof of Proposition~\ref{prop: weight}, we observe that the localization of non-scattering states is influenced by their low-frequency component, which corresponds to the $0$ threshold energy. 
\begin{theorem}\label{thm2}If Assumptions~\ref{aspV2} and~\ref{aspV3} are satisfied, then \eqref{goal: 2} holds true.
\end{theorem}

\subsubsection{The time-independent result}
When the potential is time-independent, we can write $U(t,0)=e^{-itH}$. In this case, $P_c(t_0)\equiv P_c$ is equal to the projection onto the continuous spectrum of $H$, see, for example, page $2$ of \cite{RS2004}.  
\begin{assumption}\label{aspV}$V(x)$ is real-valued and satisfies $\langle x\rangle^4 V(x)\in \s^\infty_x(\mathbb{R}^3)$.
\end{assumption}
We prove that under Assumption~\ref{aspV}, $e^{-itH}P_c$ satisfies the local decay estimate
\eq
\int_{-\infty}^\infty dt\| \langle x\rangle^{-\eta} e^{-itH}F(H\geq c)P_c\psi_0\|^2\lesssim_{\eta,c}\|\psi_0 \|^2\label{main}
\eeq
for all $\psi_0\in \mathcal L^2_x(\mathbb R^3)$, $c>0$ and any $\eta>1$.  
\begin{theorem}\label{thm1}Assume $V(x)$ satisfies Assumption \ref{aspV}. Then for any $\eta > 1$ and $c>0$, \eqref{main} holds for all initial states $\psi_0 \in \mathcal{L}^2_x(\mathbb{R}^3)$.
\end{theorem}

\begin{remark}
We believe that a similar result can be proved in four or higher space dimensions by this method.
\end{remark}

\subsection{Applications}
As an application, we obtain Strichartz estimates globally in $t$. Strichartz estimates state that
\eq
\|  U(t,0)P_c(0) f\|_{\s^q_t\s^r_x(\mathbb{R}^{n+1})}\leq C_q \| f \|_{\s^2_x(\mathbb{R}^n)}\label{Stri}
\eeq
for $2\leq r,q\leq \infty, \frac{n}{r}+\frac{2}{q}=\frac{n}{2}, \text{ and }(q,r,n)\neq (2,\infty,2).$
\begin{theorem}\label{thmS}If $V(x,t)$ satisfies Assumptions~\ref{aspV2} and~\ref{aspV3}, then Strichartz estimates(See \eqref{Stri}) are valid for all admissible $(q,r,n)$ with $n\geq 5$.
\end{theorem}

\subsection{Outline of the proof}
The proof of Theorem~\ref{thm} differs from that of Theorem~\ref{thm1}, though the essential ideas remain the same. In the time-independent case, we address the issue arising from the zero energy threshold by applying an energy cut-off $F(H\geq c)$ to the initial data. In contrast, in the time-dependent case, we resolve this issue by considering space dimensions $n \geq 5$. We provide an outline of the proof of Theorem~\ref{thm} here.

We take $t_0=0$ and write $P_c\equiv P_c(0)$ for simplicity. The case when $t_0\neq 0$ can be treated similarly. The proof is based on a new compactness argument, the notion of 'incoming/outgoing waves' (see Definition~\ref{def: in/out}) and the knowledge of AC (the identity $P_c(t)=\Omega_+(t)\Omega_+^*(t), t\in\mathbb R$). It consists of three steps. 

We introduce some notions before showing the steps. Let $P^\pm$ denote the projections onto the incoming/outgoing waves. We present the optimized (adjoint) M\"oller wave operators as introduced in \cite{SW20221}: for $\alpha\in (0, 1-\frac{2}{n}),\,n\geq 3$,
\begin{equation}\label{WO: def}
\Omega_\pm^*(t):=s\text{-}\lim\limits_{s\to \pm \infty} e^{isH_0}F(\frac{|x-2sp|}{|s|^\alpha}<1)U(s+t,t),\qquad \text{ on }\mathcal L^2_x(\mathbb R^n).
\end{equation}
Here, we use notations $\Omega_\pm^*(t)$ since $\Omega_\pm^*(t)$ are independent of $\alpha$. See Eq.~(2.13) of \cite{SW20221}. $\Omega_\pm^*(t)$ satisfy the intertwining property:
\begin{equation}
    \Omega_\pm^*(t)U(t,0)=e^{-itH_0}\Omega_\pm^*(0).\label{inter}
\end{equation}

First, we write $\psi(t)$ as the sum of a 'free flow' and a compact operator acting on $\psi(t)$. Decomposing $\psi(t)$ into the sum of $P^+\psi(t)$ and $P^-\psi(t)$ and using $P^\pm\Omega_\pm^*\psi(t)$ to approximate $P^\pm\psi(t)$, respectively, we obtain with $\psi(t)\equiv U(t,0)P_c(0)\psi_0$,
\begin{equation}
    \psi(t)=\psi_f(t)+\mathcal C(t)\psi(t),\label{psi: f+cpsi}
\end{equation}
where $\psi_f(t)$ is given by, by Eq.~\eqref{inter}, 
\begin{align}
    \psi_f(t):=&P^+ \Omega_+^*(t)\psi(t)+P^-\Omega_-^*(t)\psi(t)\nonumber\\
    =& P^+e^{-itH_0}\Omega_+^*(0)\psi_0+ P^-e^{-itH_0}\Omega_-^*(0)\psi_0\label{def: psif}
\end{align}
and $\mathcal C(t)$ is defined by 
\begin{equation}
    \mathcal C(t):=P^+(\mathbbm1-\Omega_+^*(t))+P^-(\mathbbm1-\Omega_-^*(t)).\label{def: Ct}
\end{equation}
We note that $\mathcal C(t), t\in \mathbb R$, are compact on $\mathcal L^2_x(\mathbb R^n)$ and satisfies, by estimates $\|P^\pm \|\leq 1$ and $\|\Omega^*_\pm(t)\|\leq 1$, 
\begin{equation}
    \sup\limits_{t\in \mathbb R} \|\mathcal C(t)\|\leq 4.\label{C: unbd}
\end{equation}
It is also worth noting that $\psi_f(t)$ behaves like a free flow and satisfies $L^2$ local decay estimate and $L^2$ local smoothing estimate.

Next, we use $P_c(t)=\Omega_+(t)\Omega_+^*(t)$ to decompose $\mathcal C(t)$ further. We define
\begin{equation}
    F_M(x,p):= F(|x|<M)  F(|p|\geq\frac{1}{M})\label{def: FMxp}
\end{equation}
and write with $P_b(t):=\mathbbm1-P_c(t)$,
\begin{equation}
    \mathcal C(t)=\mathcal C_M(t)+\mathcal C_r(t)+\mathcal C(t)P_b(t),
\end{equation}
where operators $\mathcal C_M(t)$ and $\mathcal C_r(t)$ are given by
\begin{equation}
\mathcal C_M(t):= \mathcal C(t)\Omega_+(t) F_M(x,p)\Omega_+^*(t)
\end{equation}
and 
\begin{equation}
    \mathcal C_r(t):=\mathcal C(t)\Omega_+(t)\left(\mathbbm1- F_M(x,p)\right)\Omega_+^*(t).
\end{equation}
Then Eq.~\eqref{psi: f+cpsi} is equivalent to 
\begin{equation}
    \psi(t)=\psi_f(t)+\mathcal C_M(t)\psi(t)+\mathcal C_r(t)\psi(t).\label{n: f+cpsi}
\end{equation}
We prove in Proposition~\ref{prop: 4.1} that there exists a constant $M_0>1$ such that whenever $M\geq M_0$,
\begin{equation}
    \sup\limits_{t\in \mathbb R} \| \mathcal C_r(t)\|<\frac{1}{2}.\label{cpt: est}
\end{equation}
The proof of Proposition~\ref{prop: 4.1} relies on the condition that $V$ is quasi-periodic in $t$ and a standard compactness argument. It is also worth noting that $\mathcal C(t)$ is quasi-periodic in $t$ by Assumption~\ref{aspV2} (see also Corollary~\ref{cor: bd: C}) and 
\begin{align}
 \int_0^\infty   \| F_M(x,p)\Omega_+^*(t)\psi(t)\|^2dt=& \int_0^\infty   \| F_M(x,p)e^{-itH_0}\Omega_+^*(0)\psi_0\|^2dt\nonumber\\
 \lesssim_M & \|\psi_0\|^2,
\end{align}
which together with~\eqref{C: unbd}, implies, with $\psi(t)=P_c(t)\psi(t)$,
\begin{equation}
    \| \mathcal C_M(t) \psi(t)\|=\| \mathcal C_M(t)P_c(t) \psi(t)\|\leq \sup\limits_{s\in \mathbb R} \| \mathcal C(s)\| \| F_M(x,p)\Omega_+^*(t)\psi(t)\|\in L^2_t[0,\infty).\label{intro: est: CM}
\end{equation}

In the third step, by employing~\eqref{cpt: est}, we find that $(\mathbbm1-\mathcal C_r(t))^{-1}$ is bounded on $\mathcal L^2_x(\mathbb R^n)$. Moving $\mathcal C_r(t)\psi(t)$ to the left-hand side of Eq.~\eqref{n: f+cpsi} and then applying $(\mathbbm1-\mathcal C_r(t))^{-1}$ to both sides, we arrive at 
\begin{equation}
    \psi(t)=(\mathbbm1-\mathcal C_r(t))^{-1}\psi_f(t)+(\mathbbm1-\mathcal C_r(t))^{-1}\mathcal C_M(t)\psi(t). \label{intro: psit decom}
\end{equation}
We prove in Proposition~\ref{BP2} that 
\begin{equation}
\|\langle x\rangle^{-\eta}(\mathbbm1-\mathcal C_r(t))^{-1}\psi_{f}(t)\|_{\mathcal L^2_{x,t}(\mathbb R^{n}\times \mathbb R^+)}\lesssim \|\psi(0)\|,\qquad n\geq 5
\end{equation}
for all $\eta>\frac{3}{2}$. This together with estimate~\eqref{intro: est: CM} and Eq.~\eqref{intro: psit decom} yields the estimate~\eqref{main}.

\subsection{Challenges of the Analysis of Time-Dependent Hamiltonians}
Time Dependent Hamiltonians cannot be treated directly by the standard spectral theory methods. The spectral properties of $H(t)$ for each fixed $t$ cannot be directly related to the solution of the the Schr\"odinger equation with Hamiltonian $H(t).$
We focus here on the problem where the interaction term is given by a localized in space potential, in 3 or higher dimensions.
We are interested in the basic questions about the behavior of the solutions: localized (bound) states, scattering states, decay estimates for the scattering states etc...

\subsubsection{Previous methods} There are classical  methods of proving dispersive estimates for linear Schr\"odinger equations: Mourre's Method and Resolvent estimates are general powerful such methods. Resolvent estimates, going back many decades ago, are based on constructing a good estimate on the resolvent of the full Hamiltonian in terms of the approximate free Hamiltonian. This approach is pretty general, but limited to localized and time-periodic interactions. Its application to N-body problems is very limited. It has no direct analog for dealing with time dependent potentials.
A more general theory is the abstract Mourre's method. It begins with the construction of a self-adjoint operator which has a positive commutator, at least when the Hamiltonian is localized near favorable points in the continuous spectrum of $H.$
This is an abstract condition. Together with extra domain assumptions, such an estimate implies Local Decay estimates, via the method of differential inequalities developed by Mourre \cite{M1979,M1981 }. A more general theory, \cite{HSS1999} allows to prove besides Local Decay, the optimal minimal and maximal velocity bounds, which are time dependent bounds.
This approach does not work at the thresholds of the Hamiltonian. It does not have a generalization to time dependent potentials in general. When $V$ is time-periodic, we use Floquet theory to reduce the problem to the one with time-independent Hamiltonian \cite{GJY2004}. We refer to such Hamiltonian as the Floquet operator. If we have a potential that is periodic in time, the corresponding Floquet Operator is time independent, with infinitely many thresholds points of the form $n\omega, n\in \mathbb{Z}$, where $\omega:=\frac{2\pi}{T}$ is the time frequency. However, in the quasi-periodic case the set of threshold points of the associated Floquet operator is dense in $\mathbb R$ and so one cannot use Mourre's Method at any point.
Our goal in this work is to introduce a new method that will cover the case of quasi-periodic time dependent potentials (in 5 or more dimensions).

\subsubsection{Connecting Quasi-Periodic Time Dependence to General Time Dependence}

 { If the Fourier transform of $V(x,t)$ in $t$ is a finite measure, then by general principles based on Wiener's theorem, $V$ can be written as a sum of almost periodic potential and a part that decays in time. The decaying part will be small eventually and the main part of the almost periodic potential is quasi-periodic, so we focus on the quasi periodic case.}


\section{Phase Space Decompositions}
In this section, we introduce the notion of incoming and outgoing waves, prove estimates for the free wave, and present properties of the evolution operator $U(t,s)$. 
\subsection{Incoming/outg     oing waves}\label{inout}
We define the incoming/outgoing wave decompositions, inspired by Mourre \cite{M1979}, based on the dilation generator $A$, given by
\eq
A:=\frac{1}{2}(x\cdot p+p\cdot x).
\eeq
\begin{definition}[Incoming/outgoing waves]\label{def: in/out} The projection on outgoing waves is defined in \cite{S2011}:
\eq
P^+:= (\tanh(\frac{A}{R})+\mathbbm1)/2\label{P+}
\eeq
with $R=100$. Here, $100$ stands a sufficiently large number. We define the projection onto incoming waves as the complement of $P^+$
\eq
P^-:=\mathbbm1-P^+.\label{P-}
\eeq
\end{definition}
In what follows, we take $R=100$. We could commute through $P^\pm$ with $|p|^a$, $|x|^a$ and $\langle x\rangle^{-1/2-\epsilon}|p|^{1/2}$ for $a\in (0,\pi R/2 )$ and $\epsilon\in (0,\frac{1}{2})$ in the sense of the lemma listed below.
\begin{lemma}\label{PPpm}For all $a\in (0, \pi R/2 )$, $\epsilon\in (0,\frac{1}{2})$, and for all $f\in H^a(\mathbb{R}^n)$, $P^\pm f$ satisfy the following:
\eq
\| |p|^aP^\pm  f\|\lesssim_{\pi/2-a/R,a/R} \|  |p|^af\|,\label{paf}
\eeq
\eq
\| |x|^aP^\pm  f\|\lesssim_{\pi/2-a/R,a/R} \|  |x|^af\|\label{Nov.20.6}
\eeq
and
\begin{equation}
    \| \langle x\rangle^{-1/2-\epsilon}|p|^{1/2}P^\pm f\|\lesssim_{\pi/2-a/R,a/R} \|\langle x\rangle^{-1/2-\epsilon}|p|^{1/2} f\|.\label{Aug.8}
\end{equation}
\end{lemma}
\begin{proof}We define a function $g$ on $\mathbb R$:
\eq
g(k):=(\tanh(k)+1)/2,\quad k\in \mathbb{R}.
\eeq
We estimate $\||p|^\alpha P^+ f\|$. $\||p|^\alpha P^- f\|$ can be treated similarly. To estimate~$\| |p|^\alpha P^+ f\|$, we find  
\begin{equation}
    \| |p|^\alpha P^+ f\|\leq \| [|p|^\alpha, P^+] f\|+\| P^+|p|^af\|.\label{cm: eq1}
\end{equation}
For the first term of the right-hand side of~\eqref{cm: eq1}, we first use the Fourier inversion theorem to express $P^+$ as
\begin{equation}
    P^+= (\tanh(\frac{A}{R})+\mathbbm1)/2=\frac{1}{\sqrt{2\pi}}\int e^{iwA/R} \hat g(w)dw,\label{P+express}
\end{equation}
where $\hat g$ denotes the Fourier transform of $g$. Plugging~\eqref{P+express} into $[|p|^a,P^+]f$ and using relation
\begin{equation}
    |p|^a e^{iwA/R}=\frac{1}{e^{aw/R}}e^{iwA/R}|p|^a,
\end{equation}
we obtain
\begin{align}
[|p|^a,P^+]f=&\frac{1}{\sqrt{2\pi}}\int \hat{g}(w) [ |p|^a, e^{iwA/R}]fdw \nonumber\\
=&\frac{1}{\sqrt{2\pi}}\int  \hat{g}(w)( \frac{1}{e^{aw/R}}-1)e^{iwA/R}|p|^afdw.\label{exp: cf}
\end{align}
In view of
\eq
\hat{g}(w)=\frac{i}{2\sqrt{2\pi}\sinh(\pi w/2)}+\frac{1}{2\sqrt{2\pi}}\delta(w)\label{exp: ghat}
\eeq
and
\begin{equation}
     \int \delta(w)  ( \frac{1}{e^{aw/R}}-1)e^{iwA/R}|p|^afdw=0 \qquad \forall\, a>0,
\end{equation}
we plug~\eqref{exp: ghat} into~\eqref{exp: cf} to obtain
\begin{equation}
\| [|p|^a,P^+]f\|\leq \frac{1}{4\pi}\int  | \frac{1}{\sinh(\pi w/2)}| | \frac{1}{e^{aw/R}}-1|\||p|^af\|dw.
\end{equation}
By estimates, for all $a\in (0, \pi R/2)$,
\begin{align}
    \int_{|w|>1}  | \frac{1}{\sinh(\pi w/2)}| | \frac{1}{e^{aw/R}}-1| dw \leq & \int_{|w|>1} \frac{2e^{-\pi w/2}}{1-e^{-\pi }} (e^{a|w|/R}+1)dw\lesssim_{\pi/2-a/R} 1
\end{align}
and
\begin{align}
     \int_{|w|\leq 1}  | \frac{1}{\sinh(\pi w/2)}| | \frac{1}{e^{aw/R}}-1| dw \leq &2\sup\limits_{|w|\leq 1} |\frac{w}{\sinh(\pi w/2)}| \sup\limits_{|w|\leq 1} |\frac{e^{-aw/R}-1}{w}|\nonumber\\
     \lesssim_{a/R} & 1, 
\end{align}
this yields 
\begin{equation}
    \|[|p|^a, P^+]f\|\lesssim_{\pi/2-a/R, a/R} \||p|^af\|\qquad \forall\, a\in (0, \pi R/2).\label{com: Paf: est}
\end{equation}
Estimate~\eqref{com: Paf: est} together with estimates~$\|P^+|p|^a f\|\leq \||p|^af\|$ and~\eqref{cm: eq1} yields~\eqref{paf} for $P^+$. Similarly, we have~\eqref{paf} for $P^-$, and proceeding as in~\eqref{cm: eq1}-\eqref{com: Paf: est} in the dual space (Fourier space), we obtain~\eqref{Nov.20.6}. By proceeding as in \eqref{cm: eq1}-\eqref{com: Paf: est} in the dual space and subsequently in the configuration space, we obtain~\eqref{Aug.8}.\end{proof}

The incoming/outgoing projections $P^\pm$ and the free flows $e^{\pm i tH_0}, t\geq 0,$ satisfy the estimate
\begin{equation}
    \|P^\pm e^{\pm itH_0} F(|p|\geq 1)F(p<2)\langle x\rangle^{-\delta}\|\lesssim_{\delta,n} \frac{1}{\langle t\rangle^\delta} ,\qquad\forall\, \delta\geq 0\text{ and }\forall\, t\geq 0. \label{est: Mourre}
\end{equation}
Estimate~\eqref{est: Mourre} follows by Mourre estimate, see~\cite{M1979,M1981, HSS1999}. By~\eqref{est: Mourre}, we obtain estimates related to $P^\pm e^{\pm itH_0}$ listed in Lemma~\ref{out/in1} below. These estimates will be the main components in proving Theorems~\ref{thm},~\ref{thm2} and~\ref{thm1}. 
\begin{lemma}\label{out/in1}When $R=100$, the operators $P^\pm e^{\pm itH_0}, t> 0,$ satisfy: given $t>0$, (recall that $n$ is the space dimension)
\begin{enumerate}
\item {\bf High Energy Estimate}
For all $\delta> 0$ and $c>0$,
\eq
\|P^{\pm}  F(|p|\geq c)e^{\pm i tH_0}\langle x\rangle^{-\delta}\|\lesssim_{n,c,\delta}  \frac{1}{\langle t\rangle^{\delta}} .\label{Sep29.1}
\eeq
\item {\bf Pointwise Smoothing Estimate} For $\delta> 0$, $l\in[0,\delta)$, $M\geq1$ and $t>0$ with $tM^2\geq 1$,
\eq
\|P^{\pm}   F(|p|\geq M)e^{\pm i tH_0}|p|^l\langle x\rangle^{-\delta}\|\lesssim_{n,\delta-l}  \frac{1}{M^{\delta-l}t^\delta} .\label{Sep20.2}
\eeq

\item {\bf Time Smoothing Estimate} For $\delta> 2$, $l=1,2,$
\eq
\int_0^1  t^2\|P^{\pm}   F(|p|\geq 1)e^{\pm i tH_0}|p|^l\langle x\rangle^{-\delta}\|dt\lesssim_{n,\delta} 1 .\label{Oct.1}
\eeq

\item {\bf Near Threshold Estimate} For all $\delta>0$ and all $\epsilon\in (0,1/2)$, $\alpha\in [0, \pi R/2)$,
\eq
\| \frac{|p|^\alpha}{\langle p\rangle^\alpha}P^{\pm}e^{\pm i tH_0} \langle x\rangle^{-\delta}\|\lesssim_{n,\epsilon,\alpha,R}  \frac{1}{\langle t\rangle^{\min\{\frac{\delta}{2}+\epsilon, (1/2-\epsilon)(\alpha+\min\{\frac{n}{2}-\epsilon,\delta\})\}}} ,\label{Sep20.1}
\eeq
and in particular, when $n\geq 5,$ $\delta>2$ and $\epsilon \in (0, \frac{n}{4}-1)$, we have  
\eq
\| P^{\pm}  e^{\pm i tH_0} \langle x\rangle^{-\delta}\|\in \s^1_t(\mathbb{R}^+).
\eeq

\item {\bf Pointwise Local Smoothing Estimate} For all $\epsilon\in (0, \frac{1}{2}), \delta\in [0, 2], n\geq 5, l=1,2,\, j=1,\cdots,n$ and $t\geq 1$,
\begin{equation}
    \|  \frac{1}{\langle x\rangle^{2-\delta}}P^\pm e^{\pm itH_0}p_j^l \langle x\rangle^{-(\frac{n}{2}+l+2-\delta)}\|\lesssim_{n,\epsilon} \frac{1}{\langle t\rangle^{\frac{n}{4}+\frac{l+2-\delta}{2}-\epsilon}}.\label{est: Pbweight}
\end{equation}
\end{enumerate}

\end{lemma}
\begin{proof} See Appendix~\ref{sec: app} for its proof.\end{proof}

\subsection{Quasi-periodic evolution operators}
We present several properties related to the quasi-periodicity of the evolution operators $U(t,s)$, for $t,s \in \mathbb{R}$, as well as the projections $P_c(t)$ and the operators $\mathcal{C}(t)$.
\begin{lemma}\label{lem: quasi}If Assumption~\ref{aspV2} holds true, then for each $t\in \mathbb R$, $U(t+s,s)$ is an evolution operator that is quasi-periodic in $s$, with the same type as $V$.
\end{lemma}
\begin{proof}Let 
\begin{equation}
    \mathscr K_t(V):=e^{itH_0}Ve^{-itH_0},\qquad t\geq 0
\end{equation}
and for a family of operators on $\mathcal L_x^2(\mathbb R^n)$, $\{A_j\}_{j=1}^{j=J}, J\in \mathbb N^+$, we define 
\begin{equation}
     \prod\limits_{j=1}^0 A_j=\mathbbm 1\qquad \text{ and }\qquad \prod\limits_{j=1}^J A_j=A_1 \cdots A_J.
\end{equation}
By the Duhamel's principle and iterating it for infinitely many times, we obtain 
\begin{equation}
    U(t+s,s)=\sum\limits_{j=0}^\infty I_j(t,s),
\end{equation}
where $I_0(t,s)=\mathbbm 1$ and $I_j(t,s), j=1,\cdots,$ are given by
\begin{equation}
   I_j(t,s):=(-i)\int_0^t \int_0^{t_1}\cdots \int_0^{t_{j-1}} e^{-itH_0}\prod\limits_{k=1}^{j} \mathscr K_{t_k}(V(x,t_k+s))  dt_j\cdots dt_1. 
\end{equation}
By Definition~\ref{def: quasif}, we conclude that operators $I_j(t,s), j=1,\cdots,$ are quasi-periodic in $s$ with the same type as $V$. Hence, $U(t+s,s)$ is quasi-periodic in $s$ with the same type as $V$. \end{proof}
Recall that 
\begin{equation}
   \Omega_\pm^*(t):= s\text{-}\lim\limits_{s\to \pm \infty} e^{isH_0}F(\frac{|x-2sp|}{|s|^\alpha}<1)U(s+t,t),\qquad \text{ on }\mathcal L^2_x(\mathbb R^n),
\end{equation}
\begin{equation}
P_c(t):=s\text{-}\lim\limits_{s\to \pm \infty} U(t,t+s)F(\frac{|x-2sp|}{|s|^\alpha}<1)U(s+t,t),\qquad \text{ on }\mathcal L^2_x(\mathbb R^n),
\end{equation}
\begin{equation}
    P_b(t):=\mathbbm1-P_c(t)
\end{equation}
and
\begin{equation}
    \mathcal C(t):=P^+(1-\Omega_+^*(t))+P^-(1-\Omega_-^*(t)).\label{mcal Ct}
\end{equation}
\begin{corollary}\label{cor: bd: C} If Assumption~\ref{aspV2} holds true, then for each $t\in \mathbb R$, $\Omega_\pm^*(t), \Omega_\pm(t), P_c(t), P_b(t)$ and $\mathcal C(t)$ are bounded operators on $\mathcal L^2_x(\mathbb R^n)$ uniformly in $t$ and quasi-periodic in $t$ with the same type as $V$. 
    
\end{corollary}
\begin{proof}The boundedness of $P_c(t)$ on $\mathcal L^2_x(\mathbb R^n)$ is proved in Theorem 2.1 of~\cite{SW20221}. The boundedness of $P_b(t)$ on $\mathcal L^2_x(\mathbb R^n)$ follows from the boundedness of $P_c(t)$, since $P_b(t) = \mathbbm{1} - P_c(t)$. The boundedness of $\mathcal C(t)$ on $\mathcal L^2_x(\mathbb R^n)$ follows from the existence of $\Omega^*_\pm(t)$ and the estimate $\|P^\pm\|\leq 1$. $P_c(t),$ $P_b(t)$ and $\Omega_\pm^*(t)$ are quasi-periodic in $t$ with the same type as $V$ by Lemma~\ref{lem: quasi}. Therefore, so is $\mathcal C(t)$.\end{proof}

\section{The time-independent problem}
In this section, we prove Theorem~\ref{thm1}. The proof is based on the properties of operator $\mathcal C(t)$ defined in Eq.~\eqref{def: Ct}. When $V$ is time-independent, $\mathcal C\equiv \mathcal C(t)$,  is time-independent. That is, 
\begin{equation}
   \mathcal C=P^+(\mathbbm1-\Omega_+^*)+P^-(\mathbbm1-\Omega_-^*).
\end{equation}
$P_c\equiv P_c(t)$ is also time-independent and we can express $U(t+s,s)=e^{-itH}$ for all $s,t\in \mathbb R$. We first derive a representation of $\psi(t)\equiv e^{-itH}F(H\geq c)P_c\psi(0)$ analogous to Eq.~\eqref{psi: f+cpsi}. Let $F_{c,1}(z):=F(z\geq\frac{c}{10})$ and $F_{c,2}(z):=F(z\geq\frac{c}{100})$. It is worth noting that Eq.~\eqref{psi: f+cpsi} holds true with $\psi(0)=F(H\geq c)\psi(0)$ when $V$ is time-independent. Using
\begin{equation}
    F(H\geq c)=F_{c,1}(H)F(H\geq c)
\end{equation}
$\psi(0)=F(H\geq c)\psi(0)$ and Eq.~\eqref{psi: f+cpsi}, we obtain 
\begin{equation}
    \psi(t)=\psi_f(t)+\mathcal C_c\psi(t),\label{1: time-in}
\end{equation}
where operator $\mathcal C_c$ is given by 
\begin{equation}
    \mathcal C_c:=\mathcal C F_{c,1}(H)
\end{equation}
and $\psi_f(t)$ reads, according to Eq.~\eqref{def: psif}, 
\begin{align}
    \psi_f(t)=P^+e^{-itH_0}\Omega_+^*\psi(0)+P^-e^{-itH_0}\Omega_-^*\psi(0).
\end{align}

In what follows in this section, we prove Theorem~\ref{thm1} by showing the compactness of $\mathcal C_c$ and employing a new compactness argument based on AC. We also assume in this section that $V$ satisfies Assumption~\ref{aspV} and $n=3$ when the context is clear.
\subsection{Compactness of $\mathcal C_c$}
We prove the compactness of $\mathcal C_c$ in this subsection.
\begin{proposition}\label{maincpt}$\mathcal C_c$ is a compact operator on $\s^2_x(\mathbb{R}^3)$ for all $c>0$.
\end{proposition}
The proof of Proposition~\ref{maincpt} requires the well-known result Lemma~\ref{lem: cpt F}
\begin{lemma}\label{lem: cpt F} $F_{c,1}(H)-F_{c,1}(H_0)$ is compact on $\mathcal L^2_x(\mathbb R^3)$ for all $c>0$. 
\end{lemma}
\begin{proof} See Appendix~\ref{sec: app} for its proof.
    
\end{proof}

\begin{proposition}\label{lem: key}For all $t\geq 1$ the estimates
\begin{equation}
    \| P^\pm F_{c,1}(H_0)e^{\pm i tH_0}\Omega_\pm^* \langle x\rangle^{-2}\|\lesssim_c \frac{1}{\langle t\rangle^{2}}
\end{equation}
are valid.
\end{proposition}
\begin{proof}Let $\mathcal O^\pm(t)\equiv P^\pm F_{c,1}(H_0)e^{\pm itH_0}\Omega_\pm^*\langle x\rangle^{-2}$. We estimate $\mathcal O^+(t)$. The estimate for $\mathcal O^-(t)$ can be derived similarly. By the Duhamel's principle, $\mathcal O^+(t)$ reads 
\begin{equation}
    \mathcal O^+(t)=\mathcal O^+_1(t)+\mathcal O^+_2(t),\label{decom: Ot}
\end{equation}
where the operators $\mathcal O^+_j(t), j=1,2,$ are given by 
\begin{equation}
    \mathcal O^+_1(t):=P^+F_{c,1}(H_0) e^{itH_0}  \langle x\rangle^{-2}
\end{equation}
and 
\begin{equation}
    \mathcal O^+_2(t):=(-i)\int_0^\infty P^+F_{c,1}(H_0) e^{i(t+s)H_0}Ve^{-isH}  \langle x\rangle^{-2}ds.
\end{equation}
By estimate~\eqref{Sep29.1}, $\mathcal O_1^+(t)$ and $\mathcal O_2^+(t)$ satisfy 
\begin{equation}
    \|\mathcal O_1^+(t)\|\lesssim_c \frac{1}{\langle t\rangle^2}
\end{equation}
and
\begin{align}
     \|\mathcal O_2^+(t)\|\leq& \int_0^\infty  \|P^+F_{c,1}(H_0)e^{i(t+s)H_0}\langle x\rangle^{-3}\|\|\langle x\rangle^3 V\|_{2\to 2} \|e^{-isH}\| ds\nonumber\\
     \lesssim_c & \int_0^\infty \frac{1}{\langle t+s\rangle^3} \| \langle x\rangle^3 V(x)\|_{\mathcal L^\infty_x(\mathbb R^3)}ds \nonumber\\
     \lesssim_c& \frac{1}{\langle t\rangle^2}\| \langle x\rangle^3 V(x)\|_{\mathcal L^\infty_x(\mathbb R^3)}, 
\end{align}
where we used the unitarity of $e^{-isH}$ and estimate $\|\langle x\rangle^3 V\|_{2\to 2}\leq  \| \langle x\rangle^3 V(x)\|_{\mathcal L^\infty_x(\mathbb R^3)}$. These together with Eq.~\eqref{decom: Ot} imply 
\begin{equation}
    \| \mathcal O^+(t)\|\lesssim_c \frac{1}{\langle t\rangle^2}.
\end{equation}
Similarly, we have 
\begin{equation}
    \| \mathcal O^-(t)\|\lesssim_c \frac{1}{\langle t\rangle^2}.
\end{equation}

\end{proof}
and 
\begin{proposition}\label{prop: key1} The operators $P^\pm e^{\pm itH_0}\Omega_\pm^*V(x)e^{\mp itH_0}, t>0$, are compact on $\mathcal L^2_x(\mathbb R^3)$.
\end{proposition}
\begin{proof}  We prove the compactness of $\mathcal O(t)\equiv P^+e^{itH_0}\Omega_+^* V(x)e^{-itH_0}$. $P^-e^{-itH_0}\Omega_-^* V(x)e^{itH_0}$ can be treated similarly. We decompose $P^+e^{itH_0}\Omega_+^* V(x)e^{-itH_0}$ into the sum of the low frequency part and the high frequency part
 \begin{equation}
    \mathcal O(t)= \mathcal O_{l}(t)+ \mathcal O_{h}(t),
 \end{equation}
 where $O_l(t)$ and $O_h(t)$ are given by 
\begin{equation}
     \mathcal O_{l}(t):=P^+ F(|p|< \frac{t^2+1}{ t})e^{itH_0}\Omega_+^* V(x)e^{-itH_0}
 \end{equation}
and
\begin{equation}
    \mathcal O_{h}(t):=P^+ F(|p|\geq \frac{t^2+1}{ t})e^{itH_0}\Omega_+^* V(x)e^{-itH_0}.
\end{equation}
By the intertwining property
\begin{equation}
    F(|p|< \frac{t^2+1}{ t})\Omega^*_\pm =\Omega^*_\pm F(\sqrt{HP_c}< \frac{t^2+1}{ t}),
\end{equation}
 $\mathcal O_l(t)$ reads 
 \begin{equation}
     \mathcal O_{l}(t)=P^+ e^{itH_0}\Omega_+^* F(\sqrt{HP_c}< \frac{t^2+1}{ t})V(x)e^{-itH_0}.\label{Ol: id1}
 \end{equation}
Since $F(\sqrt{HP_c}< \frac{t^2+1}{ t}) \langle x\rangle^{-1}$ is compact, due to Assumption~\ref{aspV}, $F(\sqrt{HP_c}< \frac{t^2+1}{ t^{2}}) V$ is compact and therefore \eqref{Ol: id1} implies the compactness of $\mathcal O_l(t)$ on $\mathcal L^2_x(\mathbb R^n)$.

 Next, we prove the compactness of $\mathcal O_h(t)$. It requires estimates
 \eq
 \|P^{\pm}   F(|p|\geq \frac{t^2+1}{ t}) e^{\pm i (t+s)H_0}|p|\langle x\rangle^{-2}\|\lesssim_{n}   \frac{t}{(t^2+1)(t+s)^2},\qquad \forall s\geq 0,\label{vari}
 \eeq
 which follow from~\eqref{Sep20.2}. To prove the compactness of $\mathcal O_h(t)$, by the Duhamel principle, $\mathcal O_h(t)$ reads 
\begin{equation}
    \mathcal O_h(t)=\mathcal O_{h1}(t)+\mathcal O_{h2}(t),
\end{equation}
where $\mathcal O_{hj}(t)$, $j=1,2,$ are given by 
\begin{equation}
    \mathcal O_{h1}(t):=P^+ e^{itH_0}  F(|p|\geq \frac{t^2+1}{t}) V
\end{equation}
and
\begin{equation}
    \mathcal O_{h2}(t):=(-i)\int_0^\infty P^+ e^{i(t+s)H_0}  F(|p|\geq \frac{t^2+1}{t}) V e^{-isH} Vds.
\end{equation}
Since $\langle x\rangle^2 |p|^{-1} F(|p|\geq \frac{1}{8}) \langle x\rangle^{-3}$ is bounded on $\mathcal L^2_x(\mathbb R^3)$ and compact, by writing  
\begin{align}
    \mathcal O_{h1}(t)=&P^+ e^{itH_0}  F(|p|\geq \frac{t^2+1}{t})  F(|p|\geq \frac{1}{8})V\nonumber\\
     =&P^+ e^{itH_0}  F(|p|\geq \frac{t^2+1}{t}) |p|\langle x\rangle^{-2} \langle x\rangle^2 |p|^{-1}F(|p|\geq \frac{1}{8})V
\end{align}
and using estimate~\eqref{vari},  we conclude that $  \mathcal O_{h1}(t)$ is compact. Similarly, $P^+ e^{i(t+s)H_0}  F(|p|\geq\frac{t^2+1}{t}) V e^{-isH} $ is compact. By estimate~\eqref{vari}, we also obtain 
\begin{equation}
    \| P^+ e^{i(t+s)H_0}  F(|p|\geq \frac{t^2+1}{t}) V e^{-isH}  \| \in L^1_s[0,\infty).
\end{equation}
This together with the compactness of $P^+ e^{i(t+s)H_0}  F(|p|\geq \frac{t^2+1}{t}) V e^{-isH} $ yields the compactness of $\mathcal O_{h2}(t)$. The compactness of $\mathcal O_{h1}(t)$ and $\mathcal O_{h2}(t)$ implies the compactness of $\mathcal O_h(t)$. And the compactness of $\mathcal O_{h}(t)$ and $\mathcal O_{l}(t)$ implies the compactness of $\mathcal O(t)$.\end{proof}

We also need the following Corollary for the proof of Theorem~\ref{thm}.
\begin{corollary}\label{cor: cpt}If $\langle x\rangle^2 V(x,t)\in \mathcal L^\infty_{x,t}(\mathbb R^{n+1}), n\geq 1$, then $P^\pm e^{isH_0}V(x,t+s), s>0, t\in \mathbb R,$ are compact operators on $\mathcal L^2_x(\mathbb R^n)$.
    
\end{corollary}
\begin{proof}It suffices to show that $\mathcal B^\pm \equiv P^\pm e^{isH_0}\langle x\rangle^{-2}$ are compact. Since $P^{\pm}   F(|p|< \frac{s^2+1}{ s}) e^{\pm i sH_0}\langle x\rangle^{-2}$ are compact, it suffices to show the compactness of $P^{\pm}   F(|p|\geq \frac{s^2+1}{ s}) e^{\pm i sH_0}\langle x\rangle^{-2}$. For this,  according to estimate~\eqref{vari}, we have 
\begin{equation}
    \|P^{\pm}   F(|p|\geq \frac{s^2+1}{ s^{2}}) e^{\pm i sH_0}|p|\langle x\rangle^{-2}\|\lesssim_{n}   \frac{s}{(s^2+1)s^2}.
\end{equation}
This together with 
\begin{align}
   & P^{\pm}   F(|p|\geq \frac{s^2+1}{ s}) e^{\pm i sH_0}\langle x\rangle^{-2}\nonumber\\
   =&\left(P^{\pm}   F(|p|\geq \frac{s^2+1}{ s}) e^{\pm i sH_0}|p| \langle x\rangle^{-2}\right)\left( \langle x\rangle^2 |p|^{-1} F(|p|\geq \frac{s^2+1}{ 10s})\langle x\rangle^{-2}\right)
\end{align}
and the compactness of $\langle x\rangle^2 |p|^{-1} F(|p|\geq \frac{s^2+1}{ 10s})\langle x\rangle^{-2}$ yields the compactness of $\mathcal B^\pm$. 
\end{proof}

\begin{proof}[Proof of Proposition~\ref{maincpt}] Note that 
\begin{equation}
    \mathcal C_c=\mathcal C_{c,+}+\mathcal C_{c,-},\label{eq: Cdeco}
\end{equation}
where the operators $\mathcal C_{c,\pm}$ are given by 
\begin{equation}
    \mathcal C_{c,\pm}:=P^\pm (\mathbbm1-\Omega_\pm^*)F_{c,1}(H).\label{eq: Cpmdeco}
\end{equation}
We prove that $\mathcal C_{c,+}$ is compact on $\mathcal L^2_x(\mathbb R^3)$ and the compactness of $\mathcal C_{c,-}$ will follow similarly. To prove the compactness of $\mathcal C_{c,+}$, we note by intertwining, 
\begin{equation}
    \mathcal C_{c,+}=P^+F_{c,1}(H)-P^+ F_{c,1}(H_0)\Omega_+^*.
\end{equation}
By Lemma~\ref{lem: cpt F}, it suffices to show that $P^+F_{c,1}(H_0)(\mathbbm1-\Omega_+^*)$ is compact. For this, using relation $\mathbbm{1}=\Omega_+^*\Omega_+$ and the Duhamel principle to expand $\Omega_+$, we obtain for all $f\in \mathcal L^2_x(\mathbb R^3)$,
\begin{align}
   P^+F_{c,1}(H_0)-P^+ F_{c,1}(H_0)\Omega_+^*=& -i\int_0^\infty P^+F_{c,1}(H_0)e^{itH_0}\Omega_+^*V(x)e^{-itH_0}dt. 
\end{align}
By Proposition~\ref{prop: key1} and the compactness of $P^+e^{itH_0}\Omega^*_+(1-F_{c,1}(H))\langle x\rangle^{-2}$, we conclude that 
\begin{equation}
\int_0^M P^+F_{c,1}(H_0)e^{itH_0}\Omega^*_+V(x)e^{-itH_0}dt \text{ is compact for each }M > 0. 
\end{equation}
On the other hand, by Proposition~\ref{lem: key} and the unitarity of $e^{-itH_0}$, we obtain 
\begin{align}
     &\| \int_M^\infty P^+ F_{c,1}(H_0)e^{ itH_0}\Omega_+^* Ve^{-itH_0}dt\|\nonumber\\
     \leq & \int_M^\infty \| P^+ F_{c,1}(H_0) e^{ itH_0}\Omega_+^*\langle x\rangle^{-2}\|\| \langle x\rangle^2V(x)\|_{2\to 2} \|e^{-itH_0}\| dt\nonumber\\
     \lesssim_c & \| \langle x\rangle^{3}V(x)\|_{\mathcal L^\infty_x(\mathbb R^3)} \int_M^\infty  \frac{1}{\langle t\rangle^2} dt\nonumber\\
     \lesssim_c & \frac{1}{M}\| \langle x\rangle^{3}V(x)\|_{\mathcal L^\infty_x(\mathbb R^3)}\to 0,
\end{align}
as $M\to \infty.$ Hence, we conclude that $ P^+F_{c,1}(H_0)-P^+ F_{c,1}(H_0)\Omega_+^*$ is compact and therefore, $\mathcal C_{c,+}$ is compact. Similarly, $\mathcal C_{c,-}$ is compact. Thus, we conclude the compactness of $\mathcal C_c$.  \end{proof}

\subsection{Decomposition of the Operator $\mathcal C_c$} 
Recall that 
\begin{equation}
    F_M(x,p)= F(|x|< M)  F(|p|\geq\frac{1}{M})\label{def: FM3}
\end{equation}
and 
\begin{equation}
    \mathcal C_c=\left(P^+(\mathbbm1-\Omega_+^*(t))+P^-(\mathbbm1-\Omega_-^*(t))\right) F_{c,1}(H).
\end{equation}
We define 
\begin{equation}
\mathcal C_{M,c}:=  \mathcal C_c\Omega_+ F_M(x,p)\Omega_+^*\label{def: CM}
\end{equation}
and
\begin{equation}
\mathcal C_{r,c}:=\mathcal C_{c}\Omega_+\left(\mathbbm1- F_M(x,p)\right)\Omega_+^*. \label{def: Cr}
\end{equation}
 Since $F(H\geq c)=F(H\geq c)P_c$, then we have
 \begin{equation}
     \mathcal C_c=\mathcal C_{M,c}+\mathcal C_{r,c}.
 \end{equation}

In this section, we show that there exists $M\geq 1$ such that for all $\psi_0\in \mathcal L^2_x(\mathbb R^3)$,
\begin{equation}
  \|  \mathcal C_{M,c} e^{-itH}P_c\psi_0 \|_{\mathcal L^2_{x,t}(\mathbb R^{3}\times \mathbb R^+)}\lesssim_M \|\psi_0\|,\label{est: M}
\end{equation}
\begin{equation}
    \| \mathcal C_{r,c}\|<\frac{1}{2}\label{est: Cr0}
\end{equation}
and for all $\eta>1$ and $f\in \mathcal L^2_x(\mathbb R^3)$,
\begin{equation}
    \| \langle x\rangle^{-\eta} (\mathbbm1-\mathcal C_{r,c})^{-1}P^\pm e^{-itH_0} f \|_{\mathcal L^2_{x,t}(\mathbb R^3\times \mathbb R^+)} \lesssim_{\eta} \|f\|.\label{est: Cr}
\end{equation}
\begin{proposition}\label{prop: main time-indp}There exists $M_0\geq 1$ such that for all $M\geq M_0$, with the decomposition $\mathcal C_c=\mathcal C_{M,c}+\mathcal C_{r,c}$, estimates~\eqref{est: M}-\eqref{est: Cr} hold true for all $\eta>1$.

\end{proposition}
The proof of Proposition~\ref{prop: main time-indp} requires the well-known result Lemma~\ref{lem: Fwell}
\begin{lemma}\label{lem: Fwell}For all $c>0$ and $\epsilon\in (0,1/2)$, when $n\geq 3$,
\begin{equation}
    \| F_{c,1}(H)(\mathbbm1-F_{c,2}(H_0))|p|^{-1/2}\langle x\rangle^{1/2+\epsilon}\|\lesssim_c 1\label{lem3.2: est: 1}
\end{equation}
and
\begin{equation}
    \| (F_{c,1}(H)-F_{c,1}(H_0))F_{c,2}(H_0)|p|^{-1/2}\langle x\rangle^{1/2+\epsilon}\|\lesssim_c 1.\label{lem3.2: est: 2}
\end{equation}
    
\end{lemma}
\begin{proof} See Appendix~\ref{sec: app}. 
\end{proof}
\begin{lemma}The estimate~\eqref{est: M} holds true for all $M\geq 1$. 
    
\end{lemma}
\begin{proof}By the intertwining property and Eq.~\eqref{def: CM}, $\mathcal C_{M,c}e^{-itH}P_c\psi_0$ reads 
\begin{equation}
    \mathcal C_{M,c}e^{-itH}P_c\psi_0= \mathcal C_c\Omega_+F_M(x,p)e^{-itH_0}\Omega_+^*\psi_0. \label{Lem3.2: eq1}
\end{equation}
By $L^2$ local decay estimate and the existence of $\Omega_+^*$, we note
\begin{equation}
    \| F_M(x,p) e^{-itH_0}\Omega_+^*\psi_0\|_{\mathcal L^2_{x,t}(\mathbb R^{3+1})}\lesssim_M \|\psi_0\|.
\end{equation}
By the boundedness of $\mathcal C, F(H\geq c), \Omega_+$ and $\Omega_+^*$(see Corollary~\ref{cor: bd: C} for the boundedness of $\mathcal C$), this together with Eq.~\eqref{Lem3.2: eq1} implies
\begin{align}
\int \| \mathcal C_{M,c}e^{-itH}P_c\psi_0\|^2 dt\lesssim &\int \|F_M(x,p) e^{-itH_0}\Omega_+^*\psi_0\|^2 dt\nonumber\\
\lesssim &\int \|F(|x|< M) e^{-itH_0} F(|p|\geq \frac{1}{M})\Omega_+^*\psi_0\|^2 dt\nonumber\\
\lesssim_M & \|\psi_0\|^2,
\end{align}
that is, 
\begin{equation}
    \|\mathcal C_{M,c} e^{-itH}P_c\psi_0\|_{\mathcal L^2_{x,t}(\mathbb R^{3+1})}\lesssim_M \|\psi_0\|. 
\end{equation}
\end{proof}
\begin{proposition}\label{prop: 2} There exists $M_0\geq 1$ large enough such that the estimate~\eqref{est: Cr0} holds true for all $M\geq M_0$.
    
\end{proposition}
\begin{proof}
  Let $P_b:=1-P_c$ be the projection on the space of discrete spectrum of $H$. By Proposition~\ref{maincpt}, $\mathcal C_c$ is compact and can therefore be approximated by a finite-rank operator: with $\mathcal C_{c}P_b=\mathcal C F(H\geq c)P_b=0$,
\eq
\|\mathcal C_c-\mathcal P_{c,N}\|\leq \frac{1}{1000},\label{id: Ccappro}
\eeq
where the finite-rank operator $\mathcal P_{c,N}$, for some $ N\in \mathbb N^+$, is given by 
\begin{equation}
    \mathcal P_{c,N}:=\sum\limits_{l=1}^{N}c_{l}({\phi}_{l}(x), \cdot )_{\s^2_x(\mathbb{R}^3)}\psi_{l}(x)\label{def: PcN}
\end{equation}
with constants $ c_{l} \in \mathbb C-\{0\}, j=1,\cdots,N,$ and $\psi_l, \phi_l\in \mathcal L^2_x(\mathbb R^3)$ satisfying 
\eq
\begin{cases}
 \| \psi_{l}\|_{\s^2_x(\mathbb{R}^3)}=1, \\
   {\phi}_{l}=P_c{\phi}_{l},\\
    \| {\phi}_{l}\|_{\s^2_x(\mathbb{R}^3)}=1.
 \end{cases}
\eeq
Estimate~\eqref{id: Ccappro} implies 
\begin{equation}
    \|\left(\mathcal C_c-\mathcal P_{c,N}\right)\Omega_+ (\mathbbm1-F_M(x,p))\Omega_+^*\|\leq \|\mathcal C_c-\mathcal P_{c,N} \|\leq \frac{1}{1000}.\label{est: a}
\end{equation}
Using that
\begin{equation}\label{lim: 1-FM}
    s\text{-}\lim\limits_{M\to \infty}\mathbbm1-F_M(x,p)=0,\qquad \text{on }\mathcal L^2_x(\mathbb R^3),
\end{equation}
$N\in \mathbb N^+$ and AC ($i.e.$ $\phi_l=P_c\phi_l=\Omega_+\Omega_+^*\phi_l, l=1,\cdots,N$), we obtain that there exists $M_0>0$ such that for all $M\geq M_0$,
\begin{equation}
    \|\Omega_+(\mathbbm 1- F_M(x,p))\Omega_+^*\phi_l\|= \|(\mathbbm 1- F_M(x,p))\Omega_+^*\phi_l\|<\frac{1}{1000N|c_l|}, \qquad \forall l=1,\cdots,N. 
\end{equation}
This together with~\eqref{def: PcN} yields 
\begin{equation}
    \| \mathcal P_{c,N}\Omega_+ (\mathbbm1-F_M(x,p))\Omega_+^* \|< \frac{1}{1000}.\label{est: b}
\end{equation}
Estimates~\eqref{est: a} and~\eqref{est: b} together with Eq.~\eqref{def: Cr} yield~\eqref{est: Cr0} for all $M\geq M_0$.\end{proof}

\begin{proposition}\label{prop: CM est} For all $M\geq 1$ and $c>0$, we have
\begin{equation}
    \|\mathcal C_{M,c} |p|^{-1/2}\langle x\rangle^{1/2+\epsilon}\|\lesssim_{M} 1
\end{equation}
    holds true for all $\epsilon\in (0,1/2)$.
\end{proposition}
and 
\begin{proposition}\label{prop: 3}Let $M_0$ be as in Proposition~\ref{prop: 2}. The estimate~\eqref{est: Cr} holds true for all $M\geq M_0$.
    
\end{proposition}
The proof of Propositions~\ref{prop: CM est} and~\ref{prop: 3} requires
\begin{proposition}\label{prop: 4}
    Recall $F_{c,1}(H_0)=F(H_0\geq\frac{c}{10})$ and $F_{c,2}(H_0)=F(H_0\geq\frac{c}{100})$, $c>0$. The estimates
\begin{equation}
    \| P^\pm (\Omega^*_\pm-\mathbbm1)F_{c,1}(H) |p|^{-1/2}\langle x\rangle^{1/2+\epsilon}\|\lesssim_{c,\epsilon} 1\label{prop: 4: est}
\end{equation}
hold true for all $c\in (0,1/2)$.
\end{proposition}

The proof of Proposition~\ref{prop: 4} requires
\begin{lemma}\label{lem: Fasp} 
\begin{equation}
    \| F_{c,1}(H) (1-F_{c,2}(H_0))|p|^{-1/2}\langle x\rangle^{1/2+\epsilon}\|\lesssim_{c,\epsilon} 1 
\end{equation}
holds true for all $\epsilon\in (0,1/2)$ and $c>0$.
    
\end{lemma}
\begin{lemma}\label{lem: Fasp2} 
\begin{equation}
    \| (F_{c,1}(H)-F_{c,1}(H_0)) F_{c,2}(H_0)|p|^{-1/2}\langle x\rangle^{1/2+\epsilon}\|\lesssim_{c,\epsilon} 1 
\end{equation}
holds true for all $\epsilon\in (0,1/2)$ and $c>0$.
    
\end{lemma}
\begin{lemma}\label{lem: estO}Let 
\begin{equation}
    \mathcal O(t):=\langle x\rangle^{-3/2} e^{-itH_0}|p|^{-1/2}\langle x\rangle^{1/2+\epsilon}.\label{def: Ot}
\end{equation}
For all $\epsilon\in (0,1/2)$, $c>0$, $n\geq3$ and $f\in \mathcal L^2_x(\mathbb R^n)$, $\mathcal O(t)$ satisfies 
\begin{equation}
   \left(\int  \langle t\rangle^{-2}\|\mathcal O(t)f\|^2 dt\right)^{1/2}\lesssim_{n,\epsilon} \|f\|. \label{lem3.6: ineq}
\end{equation}

\end{lemma}
\begin{proof}To estimate~$\mathcal O\chi(|x|\leq 1)$, we note that by Hardy-Littlewood-Sobolev inequality,
\begin{equation}
    \| \mathcal O\chi(|x|\leq 1)\|\leq \|\langle x\rangle^{-3/2}|x|^{1/2} \|\||x|^{-1/2}|p|^{-1/2}\| \|\langle x\rangle^{1/2+\epsilon}\chi(|x|\leq 1)\|\lesssim 1.
\end{equation}
Next, we estimate $\mathcal O\chi(|x|>1)$. To estimate $\mathcal O\chi(|x|>1)$, we decompose $\mathcal O\chi(|x|>1)$ into $n$ pieces 
\begin{equation}
    \mathcal O\chi(|x|>1)=\sum\limits_{j=1}^n \mathcal O_j,
\end{equation}
where operators $\mathcal O_j,j=1,\cdots,n,$ are given by
\begin{equation}
\mathcal O_j:=\mathcal O\chi(|x|>1)F_j(x),\qquad j=1,\cdots,n,
\end{equation}
with $\{F_j\}_{j=1}^{j=n}$, a partition of unity, satisfying
\begin{equation}
    |x_j|\geq \frac{1}{n}|x|\qquad \text{for all }x\in \text{supp}(F_j),\, \,  j=1,\cdots,n.\label{def: partion}
\end{equation}
We estimate $\mathcal O_1$ and the estimates for $\mathcal O_j, j=2,\cdots,n,$ can be derived similarly. Using equation
\begin{align}
   & e^{-itH_0}|p|^{-1/2}\langle x\rangle^{1/2+\epsilon}F_1(x)\nonumber\\
   =&x_1e^{-itH_0}|p|^{-1/2}\langle x\rangle^{1/2+\epsilon}\frac{F_1(x)}{x_1}-2t\frac{p_1}{|p|^{1/2}}e^{-itH_0}\langle x\rangle^{1/2+\epsilon}\frac{F_1(x)}{x_1}\nonumber\\
    &+[|p|^{-1/2},x_1]e^{-itH_0}\langle x\rangle^{1/2+\epsilon}\frac{F_1(x)}{x_1},
\end{align}
we decompose $\mathcal O_1$ into three parts 
\begin{equation}
    \mathcal O_1=\sum\limits_{j=1}^3 \mathcal O_{1j},
\end{equation}
where operators $\mathcal O_{1j}, j=1,2,3,$ are given by 
\begin{equation}
    \mathcal O_{11}:=\langle x\rangle^{-3/2}x_1 e^{-itH_0}|p|^{-1/2}\langle x\rangle^{1/2+\epsilon}\frac{F_1(x)}{x_1}\chi(|x|>1),
\end{equation}
\begin{equation}
    \mathcal O_{12}:=-2t \langle x\rangle^{-3/2}\frac{p_1}{|p|^{1/2}}e^{-itH_0}\langle x\rangle^{1/2+\epsilon}\frac{F_1(x)}{x_1}\chi(|x|>1)
\end{equation}
and
\begin{equation}
   \mathcal O_{13}:=\langle x\rangle^{-3/2}[|p|^{-1/2},x_1]e^{-itH_0}\langle x\rangle^{1/2+\epsilon}\frac{F_1(x)}{x_1}\chi(|x|>1),
\end{equation}
respectively. By Hardy-Littlewood-Sobolev inequality and the unitarity of $e^{-itH_0}$, $\mathcal O_{1j}, j=1,3,$ satisfy 
\begin{equation}
    \|\mathcal O_{11}\|\leq \sup\limits_{t\in \mathbb R}\| \langle x\rangle^{-3/2}x_1|x|^{1/2}\|\| |x|^{-1/2}|p|^{-1/2}\|\|e^{-itH_0}\|\| \langle x\rangle^{1/2+\epsilon}\frac{F_1(x)}{x_1}\chi(|x|>1)\|\lesssim_n 1\label{est: O11}
\end{equation}
and with $n\geq3,$
\begin{align}
     \|\mathcal O_{13}\|\leq&\sup\limits_{t\in \mathbb R} \| \langle x\rangle^{-3/2}|x|^{1+\epsilon}\|\||x|^{-1-\epsilon}|p|^{-1-\epsilon}\|\||p|^{1+\epsilon}[|p|^{-1/2},x_1]|p|^{1/2-\epsilon}\|\nonumber\\
     &\times\|e^{-itH_0}\|\| |p|^{-1/2+\epsilon}|x|^{-1/2+\epsilon}\|\| \langle x\rangle \frac{F_1(x)}{x_1} \chi(|x|>1)\|\nonumber\\
     \lesssim_{\epsilon,n} & 1.\label{est: O13}
\end{align}
By $L^2$ local smoothing estimate, $\mathcal O_{1,2}$ satisfies 
\begin{align}
     \int \langle t\rangle^{-2}\|\mathcal O_{12}f\|^2dt\lesssim & \int \|\langle x\rangle^{-3/2} |p|^{1/2}e^{-itH_0}\frac{p_1}{|p|}\langle x\rangle^{1/2+\epsilon}\frac{F_1(x)}{x_1}\chi(|x|>1)f\|^2dt\nonumber\\
     \lesssim_\epsilon & \| \frac{p_1}{|p|}\langle x\rangle^{1/2+\epsilon}\frac{F_1(x)}{x_1}\chi(|x|>1)f\|^2\nonumber\\
     \lesssim_{n,\epsilon}& \|f\|^2
\end{align}
for all $f\in \mathcal L^2_x(\mathbb R^n)$. This together with estimates~\eqref{est: O11} and~\eqref{est: O13} yields 
\begin{equation}
    \left(\int \langle t\rangle^{-2}\|\mathcal O_{1}f\|^2dt \right)^{1/2}\lesssim_{\epsilon,n}\|f\|,\qquad f\in \mathcal L^2_x(\mathbb R^n).
\end{equation}
Similarly, we have
\begin{equation}
    \left(\int \langle t\rangle^{-2}\|\mathcal O_{j}f\|^2dt \right)^{1/2}\lesssim_{\epsilon,n}\|f\|,\qquad f\in \mathcal L^2_x(\mathbb R^n),\, j=2,3.
\end{equation}
Hence, we conclude~\eqref{lem3.6: ineq}. \end{proof}
and
\begin{lemma}\label{lem: 3.6}Recall that $F_{c,1}(H_0)=F(H_0\geq\frac{c}{10})$ and $F_{c,2}(H_0)=F(H_0\geq\frac{c}{100})$ for $c>0$. The estimates
\begin{equation}
    \| P^\pm F_{c,1}(H_0)(\Omega^*_\pm-\mathbbm1)F_{c,2}(H_0) |p|^{-1/2}\langle x\rangle^{1/2+\epsilon}\|\lesssim_{c,\epsilon} 1\label{claim: freeest}
\end{equation}
hold true for all $c>0$ and $\epsilon\in (0,\frac{1}{2})$.
\end{lemma}
\begin{proof} We write $F_{c,1}=F_{c,1}(H_0)$ and $F_{c,2}=F_{c,2}(H_0)$ in this proof if the context is clear. We estimate~$\mathcal B^+\equiv P^+F_{c,1}(\Omega_+^*-\mathbbm 1)F_{c,2} |p|^{-1/2}\langle x\rangle^{1/2+\epsilon}$. The case of $\mathcal B^-\equiv P^-F_{c,1}(\Omega_-^*-\mathbbm 1)F_{c,2} |p|^{-1/2}\langle x\rangle^{1/2+\epsilon}$ can be treated similarly. For this, by the Duhamel principle, $\mathcal B^+$ reads
\begin{equation}\label{def: B+}
    \mathcal B^+=(-i)\int_0^\infty P^+F_{c,1} e^{itH_0}\Omega_+^*V e^{-itH_0}F_{c,2} |p|^{-1/2}\langle x\rangle^{1/2+\epsilon} dt.
\end{equation}
By Proposition~\ref{lem: key}, this implies for all $f\in \mathcal L^2_x(\mathbb R^3)$,
\begin{align}
    \|\mathcal B^+f\|\leq & \int_0^\infty \| P^+ F_{c,1}e^{itH_0}\Omega_+^*\langle x\rangle^{-2}\|\|\langle x\rangle^{7/2} V\|_{2\to 2} \| \langle x\rangle^{-3/2}F_{c,2}\langle x\rangle^{3/2}\|\|\mathcal O(t)f\|dt\nonumber\\
    \lesssim_c & \int_0^\infty \frac{1}{\langle t\rangle^{2}} \|\langle x\rangle^{7/2} V(x)\|_{\mathcal L^\infty_x(\mathbb R^3)} \|\mathcal O(t)f\|dt,\label{p3.8 B+f}
\end{align}
where operator $\mathcal O(t)$ is defined in Eq.~\eqref{def: Ot}. Applying Lemma~\ref{lem: estO} and Cauchy-Schwarz inequality to~\eqref{p3.8 B+f}, we arrive at 
\begin{align}
    \|\mathcal B^+f\|\lesssim_c&\int_0^\infty \frac{1}{\langle t\rangle} \|\langle x\rangle^{7/2} V(x)\|_{\mathcal L^\infty_x(\mathbb R^3)} \left(\frac{1}{\langle t\rangle}\|\mathcal O(t)f\|\right)dt \nonumber\\
    \lesssim_c & \|\langle x\rangle^{7/2} V(x)\|_{\mathcal L^\infty_x(\mathbb R^3)} \left(\int \frac{1}{\langle t\rangle^2} \|\mathcal O(t)f\|^2dt\right)^{1/2}\nonumber\\
    \lesssim_c &  \|\langle x\rangle^{7/2} V(x)\|_{\mathcal L^\infty_x(\mathbb R^3)}\|f\|. 
\end{align}

\end{proof}
Next, we prove Proposition~\ref{prop: 4}.
\begin{proof}[Proof of Proposition~\ref{prop: 4}] Let 
\begin{equation}
    \mathcal B^\pm:=P^\pm (\Omega_\pm^*-\mathbbm1)F_{c,1}(H) |p|^{-1/2}\langle x\rangle^{1/2+\epsilon}.
\end{equation}
We estimate $\mathcal B^+$ and the estimate for $\mathcal B^-$ can be derived similarly. We decompose $\mathcal B^+$ into two parts 
\begin{equation}
    \mathcal B^+=\mathcal B^+_1+\mathcal B^+_2,
\end{equation}
where the operators $\mathcal B^+_{j}, j=1,2,$ are given by 
\begin{equation}
    \mathcal B^+_1:=P^+(\Omega_+^*-\mathbbm1)F_{c,1}(H)(1-F_{c,2}(H_0))|p|^{-1/2}\langle x\rangle^{1/2+\epsilon}
\end{equation}
and
\begin{equation}
    \mathcal B^+_2:=P^+(\Omega_+^*-\mathbbm1)F_{c,1}(H)F_{c,2}(H_0)|p|^{-1/2}\langle x\rangle^{1/2+\epsilon},
\end{equation}
respectively. By Lemma~\ref{lem: Fasp}, we obtain 
\begin{equation}
    \|\mathcal B^+_1\|\leq 2\| F_{c,1}(H)(\mathbbm1-F_{c,2}(H_0))|p|^{-1/2}\langle x\rangle^{1/2+\epsilon}\|\lesssim_c 1.\label{est: B1: time-in}
\end{equation}
Then it suffices to estimate $\mathcal B^*_2$. For this, by intertwining, we write $\mathcal B_2^+$ as 
\begin{align}
    \mathcal B_2^+=&\mathcal B_{21}^++\mathcal B^+_{22},
\end{align}
where $\mathcal B^+_{2j}, j=1,2,$ are given by 
\begin{equation}
   \mathcal B^+_{21}:= P^+F^+_{c,1}(H_0)(\Omega_+^*-\mathbbm1)F_{c,2}(H_0)|p|^{-1/2}\langle x\rangle^{1/2+\epsilon}
\end{equation}
and 
\begin{equation}
    \mathcal B^+_{22}:= P^+ (F_{c,1}(H_0)-F_{c,1}(H))F_{c,2}(H_0)|p|^{-1/2}\langle x\rangle^{1/2+\epsilon}. 
\end{equation}
By Lemma~\ref{lem: 3.6}, we have 
\begin{equation}
    \|\mathcal B_{21}^+\|\lesssim_{c,\epsilon} 1.\label{est: B21: time-in}
\end{equation}
By Lemma~\ref{lem: Fasp2}, we have 
\begin{equation}
    \|\mathcal B_{22}^+\|\lesssim_{c,\epsilon} 1.\label{est: B22: time-in}
\end{equation}
Estimates~\eqref{est: B21: time-in} and~\eqref{est: B22: time-in} yield 
\begin{equation}
    \|\mathcal B_2^+\|\lesssim_{c,\epsilon} 1.
\end{equation}
This together with estimate~\eqref{est: B1: time-in} yields 
\begin{equation}
    \|\mathcal B^+\|\lesssim_{c,\epsilon} 1.
\end{equation}
Similarly, we have 
\begin{equation}
    \| \mathcal B^-\|\lesssim_{c,\epsilon} 1. 
\end{equation}

\end{proof}

Now we prove Propositions~\ref{prop: CM est} and~\ref{prop: 3}.
\begin{proof}[Proof of Proposition~\ref{prop: CM est}] Recall that 
\begin{equation}
    \mathcal C_{M,c}=\mathcal C_c\Omega_+ F_M(x,p)\Omega_+^*.
\end{equation}
To estimate~$\mathcal C_{M,c}$, by the Duhamel's principle, $\mathcal C_{M,c}$ reads 
\begin{equation}
    \mathcal C_{M,c}=\mathcal C_{M,c1}+\mathcal C_{M,c2},\label{dec: CMc0}
\end{equation}
where 
\begin{equation}
    \mathcal C_{M,c1}:=\mathcal C_c\Omega_+ F_M(x,p)|p|^{-1/2}\langle x\rangle^{1/2+\epsilon}
\end{equation}
and
\begin{equation}
    \mathcal C_{M,c2}:=(-i)\int_0^\infty \mathcal C_c\Omega_+ F_M(x,p) e^{itH_0}\Omega_+^* V(x)e^{-itH_0}|p|^{-1/2}\langle x\rangle^{1/2+\epsilon}.\label{id: CMc2}
\end{equation}
By the definition of $F_M(x,p)$ (see Eq.~\eqref{def: FM3}) and the boundedness of $\Omega_+$ and $\mathcal C$ (see Corollary~\ref{cor: bd: C}), $\mathcal C_{M,c1}$ satisfies 
\begin{equation}
    \|\mathcal C_{M,c1}\|\lesssim \| F_{M}(x,p)|p|^{-1/2}\langle x\rangle^{1/2+\epsilon}\| \lesssim_{M} 1.\label{est: mCMc1}
\end{equation}
Next, we estimate $\mathcal C_{M,c2}$. To estimate~$\mathcal C_{M,c2}$, we note 
\begin{equation}
    \|F_M(x,p)e^{itH_0}\Omega_+^*\langle x\rangle^{-2} \|\lesssim_{M} \frac{1}{\langle t\rangle^2}.\label{est: OmegaM}
\end{equation}
Indeed we have
\begin{equation}
    \|F_{M}(x,p)  e^{itH_0}\langle x\rangle^{-3}\|\lesssim_{M} \frac{1}{\langle t\rangle^3}\qquad \forall t\in \mathbb R.
\end{equation}
Using this estimate and the unitarity of $e^{-isH}$, we obtain estimate 
\begin{align}
    \|\int_0^\infty  F_M(x,p)e^{i(t+s)H_0}Ve^{-isH} ds\|\lesssim_{M} &\int_0^\infty \frac{1}{\langle t+s\rangle^{3}}\|\langle x\rangle^3V\|_{2\to 2}\|e^{-isH}\|ds\nonumber\\
    \lesssim_{M}& \frac{1}{\langle t\rangle^2}\|\langle x\rangle^3 V(x)\|_{\mathcal L^\infty_x(\mathbb R^3)} 
\end{align}
for all $t\geq 0$. This together with 
\begin{equation}
    \|F_{M}(x,p) e^{itH_0}\langle x\rangle^{-2}\|\lesssim_{M} \frac{1}{\langle t\rangle^2}\qquad \forall t\in \mathbb R
\end{equation}
and the Duhamel's principle yields~\eqref{est: OmegaM}.\par \eqref{est: OmegaM} together with the boundedness of $\mathcal C$ (see Corollary~\ref{cor: bd: C}) and $\Omega_+$ and Eq.~\eqref{id: CMc2} yields 
\begin{align}
    \| \mathcal C_{M,c2}\|\leq& \int_0^\infty \| F_M(x,p)e^{itH_0}\Omega_+^*\langle x\rangle^{-2}\|\|\langle x\rangle^{7/2} V\|_{2\to 2} \| \mathcal O(t)\|  dt\nonumber\\
    \lesssim_{M} &\int_0^\infty \frac{1}{\langle t\rangle^2} \|\langle x\rangle^{7/2} V(x)\|_{\mathcal L^\infty_x(\mathbb R^3)} \|\mathcal O(t)\|dt.\label{est: CMc2}
\end{align}
Here, recall that 
\begin{equation}
   \mathcal O(t)= \langle x\rangle^{-3/2}e^{-itH_0}|p|^{-1/2}\langle x\rangle^{1/2+\epsilon}.
\end{equation}
Applying Lemma~\ref{lem: estO} and Cauchy-Schwarz inequality to estimate~\eqref{est: CMc2}, we arrive at 
\begin{equation}
    \|\mathcal C_{M,c2}\|\lesssim_{M,\epsilon} \|\langle x\rangle^{7/2} V(x)\|_{\mathcal L^\infty_x(\mathbb R^3)}.
\end{equation}
This together with estimate~\eqref{est: mCMc1} and Eq.~\eqref{dec: CMc0} yields 
\begin{equation}
    \|\mathcal C_{M,c} \|\lesssim_{M,\epsilon} 1. 
\end{equation}
We complete the proof. \end{proof}

\begin{proof}[Proof of Proposition~\ref{prop: 3}] Using estimate~\eqref{prop: 4: est} and Eqs.~\eqref{eq: Cdeco} and~\eqref{eq: Cpmdeco}, we obtain
\begin{equation}
     \|\mathcal C_c |p|^{-1/2}\langle x\rangle^{1/2+\epsilon}\|\lesssim_{c,\epsilon} 1\label{est: Cabs}
\end{equation}
for all $\epsilon\in(0,1/2)$. This together with Proposition~\ref{prop: CM est} and $\mathcal C_{r,c}=\mathcal C_c-\mathcal C_{M,c}$ yields 
\begin{align}
    \| \mathcal C_{r,c}P^\pm e^{-itH_0}f \|_{\mathcal L^2_{x,t}(\mathbb R^3\times \mathbb R^+)}\leq&  \left(\int_0^\infty \| \mathcal C_{r,c} P^\pm e^{-itH_0}f\|^2 dt\right)^{1/2}\nonumber\\
    \lesssim_{M,c} & \left( \int_0^\infty \| \langle x\rangle^{-1/2-\epsilon}|p|^{1/2} P^\pm e^{-itH_0}f\|^2dt\right)^{1/2}.
\end{align}
By estimates~\eqref{Nov.20.6} and~\eqref{Aug.8}, this yields 
\begin{equation}
    \| \mathcal C_{r,c}P^\pm e^{-itH_0}f \|_{\mathcal L^2_{x,t}(\mathbb R^3\times \mathbb R^+)}\lesssim_{M} \| \langle x\rangle^{-1/2-\epsilon}|p|^{1/2}  e^{-itH_0}f\|_{\mathcal L^2_{x,t}(\mathbb R^3\times \mathbb R^+)},
\end{equation}
which leads to, by $L^2$ local decay estimate and $L^2$ local smoothing estimate, 
\begin{equation}
     \| \mathcal C_{r,c}P^\pm e^{-itH_0}f \|_{\mathcal L^2_{x,t}(\mathbb R^3\times \mathbb R^+)}\lesssim \|f\|.
\end{equation}
Therefore, using $(\mathbbm1-\mathcal C_{r,c})^{-1}=\mathbbm 1+(\mathbbm1-\mathcal C_{r,c})^{-1}\mathcal C_{r,c}$ and subsequently
\begin{align}
    \| \langle x\rangle^{-\eta}(\mathbbm1-\mathcal C_{r,c})^{-1}P^\pm e^{-itH_0}f\|\leq & \| \langle x\rangle^{-\eta}P^\pm e^{-itH_0}f\|+\| \langle x\rangle^{-\eta}(\mathbbm1-\mathcal C_{r,c})^{-1}\mathcal C_{r,c}P^\pm e^{-itH_0}f\|\nonumber\\
    \lesssim & \| \langle x\rangle^{-\eta} e^{-itH_0}f\|+\|\mathcal C_{r,c} P^\pm e^{-itH_0}f\|,
\end{align}
we arrive at 
\begin{align}
     \| \langle x\rangle^{-\eta}(\mathbbm1-\mathcal C_{r,c})^{-1}P^\pm e^{-itH_0}f\|_{\mathcal L^2_{x,t}(\mathbb R^{3}\times \mathbb R^+)}\lesssim&\| \langle x\rangle^{-\eta} e^{-itH_0}f\|_{\mathcal L^2_{x,t}(\mathbb R^{3}\times \mathbb R^+)}+\|\mathcal C_{r,c} P^\pm e^{-itH_0}f\|_{\mathcal L^2_{x,t}(\mathbb R^{3}\times \mathbb R^+)}\nonumber\\
     \lesssim_\eta & \|f\|
\end{align}
for all $\eta>1$. \end{proof}

\subsection{Proof of Theorem \ref{thm1}}
\begin{proof}[Proof for Theorem \ref{thm1}] Using Eqs.~\eqref{1: time-in},~\eqref{def: CM} and~\eqref{def: Cr} and that $\mathcal C_cP_b=0$, we obtain 
\begin{equation}
    \psi(t)=\psi_f(t)+\mathcal C_{M,c}\psi(t)+\mathcal C_{r,c}\psi(t).\label{eq: psit:in}
\end{equation}
Using~\eqref{est: Cr0}, we have $(\mathbbm1-\mathcal C_{r,c})^{-1}$ is bounded on $\mathcal L^2_x(\mathbb R^3)$. Then moving $\mathcal C_{r,c}\psi(t)$ to the left-hand side of the Eq.~\eqref{eq: psit:in} and applying $(\mathbbm1-\mathcal C_{r,c})^{-1}$ to both sides of the Eq.~\eqref{eq: psit:in}, we arrive at 
\begin{equation}
    \psi(t)=(\mathbbm1-\mathcal C_{r,c})^{-1}\psi_f(t)+(\mathbbm1-\mathcal C_{r,c})^{-1}\mathcal C_{M,c}\psi(t). 
\end{equation}
This together with estimates~\eqref{est: M} and~\eqref{est: Cr} yields~\eqref{main}.
\end{proof}

\section{Time-dependent Problems}\label{sect}
In this section, we prove Theorems~\ref{thm} and~\ref{thm2}. We take the space dimension $n \geq 5$ in this section.
\subsection{Compactness of $\mathcal C(t)$ and its decomposition}
Recall that 
\begin{equation}
    \mathcal C(t):=P^+(\mathbbm1-\Omega_+^*(t))+P^-(\mathbbm1-\Omega_-^*(t)).\label{def: mCt}
\end{equation}

We prove the compactness of $\mathcal C(t)$ in this subsection.
\begin{lemma}\label{lem: t-cpt}If $ \langle x\rangle^3 V(x,t)\in \mathcal L^\infty_{x,t}(\mathbb R^{n+1})$, then $\mathcal C(t)$ is compact on $\mathcal L^2_x(\mathbb R^n)$ for all $t\in \mathbb R$. 
    
\end{lemma}
\begin{proof}It suffices to show the compactness of $P^\pm(\mathbbm 1-\Omega_\pm^*(t))$. We prove the compactness of $P^+(\mathbbm 1-\Omega_+^*(t))$ and the compactness of $P^-(\mathbbm 1-\Omega_-^*(t))$ can be derived similarly. By the Duhamel's principle, $P^+(\mathbbm 1-\Omega_+^*(t))$ reads 
\begin{equation}
    P^+(\mathbbm 1-\Omega_+^*(t))=i\int_0^\infty  P^+e^{isH_0}V(x,s+t)U(t+s,t)ds.
\end{equation}
To prove the compactness of $P^+(\mathbbm 1-\Omega_+^*(t))$, by Corollary~\ref{cor: cpt} and Assumption~\ref{aspV2}, we obtain the compactness of $i\int_0^M P^+e^{isH_0}V(x,s+t)U(t+s,t)ds$. This together with estimate
\begin{align}
    \|\int_M^\infty P^+e^{isH_0}V(x,s+t)U(t+s,t)ds \|\leq& \int_M^\infty \| P^+e^{isH_0}\langle x\rangle^{-3}\|\|\langle x\rangle^3 V(x,t)\|_{\mathcal L^\infty_{x,t}(\mathbb R^{5+1})}ds\nonumber\\
    \lesssim & \int_M^\infty \frac{1}{\langle s\rangle^2} \|\langle x\rangle^3 V(x,t)\|_{\mathcal L^\infty_{x,t}(\mathbb R^{5+1})}ds\nonumber\\
    \to & 0\label{lem cpt: time}
\end{align}
as $M\to \infty$, yields the compactness of $P^+(\mathbbm 1-\Omega_+^*(t))$. Here, in~\eqref{lem cpt: time}, we used estimate~\eqref{Sep20.1}, with $\alpha=0$, $n\geq5$ and $\delta=3$, and the unitarity of $U(t+s,t)$. \end{proof}
Next, we prove the decomposition property of $\mathcal C(t)$. Similar to the time-independent case, we note that $\mathcal C(t)$ can be expressed as the sum of three operators
\begin{equation}
    \mathcal C(t)=\mathcal C_M(t)+\mathcal C_r(t)+\mathcal C(t)P_b(t),
\end{equation}
where the operators $\mathcal C_M(t)$ and $\mathcal C_r(t)$ are given by for $M>0$,
\begin{equation}
    \mathcal C_M(t):=\mathcal C(t)\Omega_+(t)F_{M}(x,p)\Omega_+^*(t)
\end{equation}
and
\begin{equation}
    \mathcal C_r(t):=\mathcal C(t)\Omega_+(t)(\mathbbm1-F_{M}(x,p))\Omega_+^*(t).
\end{equation}
Similarly, we use the intertwining 
\begin{equation}
    \Omega_+^*(t)U(t,0)=e^{-itH_0}\Omega_+^*(0),\qquad \text{ on }\mathcal L^2_x(\mathbb R^n)
\end{equation}
to obtain 
\begin{equation}
    \mathcal C_M(t)U(t,0)f=\mathcal C(t)\Omega_+(t)F_M(x,p)e^{-itH_0}\Omega_+^*(0)f.\label{id: CMt}
\end{equation}
We show in what follows that there exists $M_0>0$ such that whenever $M\geq M_0$, 
\begin{equation}
    \sup\limits_{t\in \mathbb R} \| \mathcal C_r(t)\|<\frac{1}{2}.\label{def: Cr bd}
\end{equation}
\begin{proposition}\label{prop: 4.1}If Assumption~\ref{aspV2} is satisfied, then there exists $M_0>0$ such that whenever $M\geq M_0$, \eqref{def: Cr bd} holds true.    
\end{proposition}
The proof of Proposition~\ref{prop: 4.1} relies on Corollary~\ref{cor: bd: C}. By Corollary~\ref{cor: bd: C}, we have $P_b(t)$, $U(T+t,t), \Omega_+(t)$ and $\Omega_+^*(t)$ are quasi-periodic in $t$ with the same type as $V$. So we could use finitely many parameters $s_j, j=1,\cdots,N$ (the same as those for $V$), which are in compact set $\bar{\mathbb T}_1\times \cdots \times \bar{\mathbb T}_N$ to express these operators: with $\vec s=(s_1,\cdots,s_N)$,
\begin{equation}
    \tilde P_b(\vec s):=P_b(t), \quad \tilde U_{\vec s}(T,0):=U(T+t,t),\quad \tilde \Omega_+(\vec s):=\Omega_+(t)\quad\text{and}\quad\tilde \Omega_+^*(\vec s):=\Omega_+^*(t)\label{def: tPb}
\end{equation}
for 
\begin{equation}
    s_j=t\, \quad\text{mod}\,\, T_j, \qquad j=1,\cdots,N.\label{vec s and t}
\end{equation}
We prove Proposition~\ref{prop: 4.1} by showing 
\begin{proposition}\label{prop: 4.2}For all $\vec s\in \bar{\mathbb T}_1\times \cdots \times \bar{\mathbb T}_N$,
\begin{equation}
    s\text{-}\lim\limits_{\vec u\to \vec s} \tilde P_b(\vec u)=\tilde P_b(\vec s). 
\end{equation}
    
\end{proposition}
The proof of Proposition~\ref{prop: 4.2} requires two lemmas below.
\begin{lemma} For each pair $\vec s, \vec u\in \bar{\mathbb T}_1\times \cdots \times \bar{\mathbb T}_N$,
\begin{equation}\label{id: dimPbs}
\text{dim}\, \tilde P_b(\vec s)=\text{dim}\, \tilde P_b(\vec u)<\infty .
\end{equation}

\end{lemma}
\begin{proof} By \cite{JL1991}, we obtain that the kernel of $\Omega_\pm^*(t)$ is equal to the space of all eigenfunctions of the Floquet operator $K$ with $s=t$. This yields $\text{Ran}(\Omega_+(t))=\text{Ran}(\Omega_-(t))$. Therefore for each $\psi\in \text{Ran}(P_b(t))$, we have $\Omega_+^*(t)\psi=\Omega_-^*(t)\psi=0$, which implies
\begin{equation}
   \psi=P^+(\mathbbm 1-\Omega_+^*(t))\psi+P^-(\mathbbm1-\Omega_-^*(t))\psi=\mathcal C(t)\psi.
\end{equation}
By Lemma~\ref{lem: t-cpt}, this implies that $\psi$ is an eigenfunction of the compact operator $\mathcal C(t)$ associated with an eigenvalue $1$. Hence, $\text{dim} (P_b(t))<\infty$ for each $t\in \mathbb R$. This together with equation
\begin{equation}
    P_b(t)=\mathbbm 1-\Omega_+(t)\Omega_+^*(t)=U(t,0)\left(\mathbbm 1-\Omega_+(0)\Omega_+^*(0)\right)U(0,t)=U(t,0)P_b(0)U(0,t)\qquad \forall t\in \mathbb R
\end{equation}
yields 
\begin{equation}
    \text{dim} (P_b(t))=\text{dim} (P_b(0))<\infty,
\end{equation}
which together with Eq.~\eqref{def: tPb} and the connection between the evolution operators and the corresponding Floquet operators implies~\eqref{id: dimPbs}. \end{proof}
\begin{lemma}\label{lem: Omegas-uF}If Assumption~\ref{aspV2} is satisfied, then 
\begin{equation}
    \lim\limits_{\vec u\to \vec s} \|(\tilde \Omega_+(\vec s)-\tilde \Omega_+(\vec u))F_M(x,p)\|=\lim\limits_{\vec u\to \vec s} \|F_M(x,p)(\tilde \Omega_+^*(\vec s)-\tilde \Omega_+^*(\vec u))\|=0.
\end{equation}
    
\end{lemma}
\begin{proof} We prove 
\begin{equation}
      \lim\limits_{\vec u\to \vec s} \|(\tilde \Omega_+(\vec s)-\tilde \Omega_+(\vec u))F_M(x,p)\|=0\label{id: Omega FMr}
\end{equation}
and 
\begin{equation}
      \lim\limits_{\vec u\to \vec s} \|F_M(x,p)(\tilde \Omega_+^*(\vec s)-\tilde \Omega_+^*(\vec u))\|=0\label{id: Omega FMl}
\end{equation}
can be proven similarly. \eqref{id: Omega FMr} follows from
\begin{equation}
     \lim\limits_{\vec u\to \vec s}  \| F_M(x,p) ( e^{iTH_0}\tilde U_{\vec u}(T,0)-e^{iTH_0}\tilde U_{\vec s}(T,0)) \|=0\label{prop4.1: lim1}
\end{equation}
and
\begin{equation}
    \lim\limits_{T\to \infty} \sup\limits_{\vec u \in \mathbb T_1\times \cdots\times \mathbb T_N} \| F_M(x,p) (\tilde \Omega_+^*(\vec u) -e^{iTH_0}\tilde U_{\vec u}(T,0)) f \|=0.\label{prop4.1: lim2}
\end{equation}
Let $\tilde V_{\vec s}(x,T):=V(x,t+T)$ whenever $\vec s$ and $t$ satisfy the condition~\eqref{vec s and t}. To prove~\eqref{prop4.1: lim1}, we use the Duhamel's principle to write 
\begin{align}
    &F_M(x,p) ( e^{iTH_0}\tilde U_{\vec u}(T,0)-e^{iTH_0}\tilde U_{\vec s}(T,0))\nonumber\\
    =&(-i)\int_0^TF_M(x,p)e^{iTH_0}\tilde U_{\vec u}(T,t)(\tilde V_{\vec s}(x,t)-\tilde V_{\vec u}(x,t))\tilde U_{\vec s}(t,0)dt.
\end{align}
This together with Condition~\eqref{asp: conti} and the unitarity of $e^{iTH_0}, \tilde U_{\vec u}(T,t)$ and $\tilde U_{\vec s}(t,0)$ yields~\eqref{prop4.1: lim1}. 

To prove~\eqref{prop4.1: lim2}, by the Duhamel's principle, we obtain 
\begin{align}
    F_M(x,p) ( e^{iTH_0}\tilde \Omega_+^*(\vec u)-e^{iTH_0}\tilde U_{\vec u}(T,0)) =&(-i)\int_T^\infty F_M(x,p)e^{itH_0} \tilde V_{\vec u}(x,t)\tilde U_{\vec u}(t,0)dt.
\end{align}
This together with estimate 
\begin{equation}
    \|F_M(x,p) e^{itH_0}\langle x\rangle^{-2}\|\lesssim_M \frac{1}{\langle t\rangle^2} 
\end{equation}
and the unitarity of $\tilde U_{\vec u}(t,0)$ yields  
\begin{align}
    \| F_M(x,p) ( e^{iTH_0}\tilde \Omega_+^*(\vec u)-e^{iTH_0}\tilde U_{\vec u}(T,0))  \|\lesssim_M& \int_T^\infty \frac{1}{\langle t\rangle^2} \|\langle x\rangle^2 V(x,w)\|_{\mathcal L^\infty_{x,w}(\mathbb R^{n+1})} dt.
\end{align}
This leads to~\eqref{prop4.1: lim2}.
\end{proof}
Next, we prove Proposition~\ref{prop: 4.2}.
\begin{proof}[Proof of Proposition~\ref{prop: 4.2}] We first claim that for all $\vec{s}\in \mathbb{R}^N$,
\eq
s\text{-}\lim\limits_{\vec{u}\to \vec{s}} \tilde P_c(\vec{u})\tilde P_c(\vec{s})=\tilde P_c(\vec{s})\quad \text{ on }\s^2_x(\mathbb{R}^n)\label{Pcs'toPcs}
\eeq
and prove the claim in the end. \eqref{Pcs'toPcs} implies
\eq
s\text{-}\lim\limits_{\vec{u}\to \vec{s}} \tilde P_b(\vec{u})\tilde P_c(\vec{s})=s\text{-}\lim\limits_{\vec{u}\to \vec{s}} (\mathbbm1-\tilde P_c(\vec{u}))\tilde P_c(\vec{s})=0\quad \text{ on }\s^2_x(\mathbb{R}^n),
\eeq
which in turn implies 
\eq
s\text{-}\lim\limits_{\vec{u}\to \vec{s}} \tilde P_b(\vec{u})-\tilde P_b(\vec{u})\tilde P_b(\vec{s})=0\quad \text{ on }\s^2_x(\mathbb{R}^n).\label{Nov.20.3}
\eeq
By~\eqref{id: dimPbs}, we arrive at 
\eq
s\text{-}\lim\limits_{\vec{u}\to \vec{s}} \tilde P_b(\vec{u})-\tilde P_b(\vec{s})=0\quad \text{ on }\s^2_x(\mathbb{R}^n).
\eeq
Now we prove the claim. Take $f\in \mathcal L^2_x(\mathbb R^n)$. Given $\epsilon>0$, by~\eqref{lim: 1-FM} and estimates $\|\tilde P_c(\vec u)\|\leq 1$ and $\|\tilde \Omega_+(\vec s)\|\leq 1$, there exists $M_1=M_1(\epsilon)>0$ such that 
\begin{equation}
    \|\tilde P_c(\vec u) \tilde \Omega_+(\vec s)(\mathbbm1-F_M(x,p))\tilde \Omega_+^*(\vec s)f\|<\frac{\epsilon}{4}
\end{equation}
and 
\begin{equation}
     \| \tilde \Omega_+(\vec s)(\mathbbm1-F_M(x,p))\tilde \Omega_+^*(\vec s)f\|<\frac{\epsilon}{4}.
\end{equation}
These yield 
\begin{equation}
    \|(\tilde P_c(\vec u)-\mathbbm1) \tilde \Omega_+(\vec s)(\mathbbm1-F_M(x,p))\tilde \Omega_+^*(\vec s)f\|<\frac{\epsilon}{2}.\label{ineq: Pu-Ps}
\end{equation}
Next, we estimate
\begin{align}
    \mathscr R_{M_1}(\vec u,\vec s)f:=& \left(\tilde P_c(\vec u)-\mathbbm1\right)\tilde P_c(\vec s)f-\left(\tilde P_c(\vec u) -\mathbbm1\right)\tilde \Omega_+(\vec s)(\mathbbm1-F_M(x,p))\tilde \Omega_+^*(\vec s)f.
\end{align}
By relation~$\tilde P_c(\vec s)=\tilde \Omega_+(\vec s)\tilde \Omega_+^*(\vec s)$, this yields 
\begin{equation}
    \mathscr R_{M_1}(\vec u,\vec s)f=\left(\tilde P_c(\vec u) -\mathbbm1\right)\tilde \Omega_+(\vec s)F_M(x,p)\tilde \Omega_+^*(\vec s)f.
\end{equation}
This together with Lemma~\ref{lem: Omegas-uF} and equation
\begin{equation}
    \left(\tilde P_c(\vec u) -\mathbbm1\right)\tilde \Omega_+(\vec u)=0
\end{equation}
yields 
\begin{equation}
    \lim\limits_{\vec u\to \vec s} \| \mathscr R_{M_1}(\vec u,\vec s)f\|=0.
\end{equation}
This together with~\eqref{ineq: Pu-Ps} yield 
\begin{equation}
    \limsup\limits_{\vec u\to \vec s}\|\left(\tilde P_c(\vec u)-\mathbbm1\right)\tilde P_c(\vec s)f\|<\epsilon,\qquad \forall \epsilon>0.
\end{equation}
Thus, we conclude the claim~\eqref{Pcs'toPcs}.\end{proof}
Now we prove Proposition~\ref{prop: 4.1}.
\begin{proof}[Proposition~\ref{prop: 4.1}] Following the proof of Proposition~\ref{prop: 2}, we note that it suffices to show that for each $f\in \mathcal L^2_x(\mathbb R^5)$, 
\begin{equation}
    \lim\limits_{M\to \infty} \sup\limits_{t\in \mathbb R}\| (\mathbbm1-F_M(x,p))\Omega_+^*(t)f \| =0
\end{equation}
holds true. By Corollary~\ref{cor: bd: C}, this is equivalent to 
\begin{equation}
     \lim\limits_{M\to \infty} \sup\limits_{\vec s\in \mathbb T_1\times \cdots \times \mathbb T_N}\| (\mathbbm1-F_M(x,p))\tilde\Omega_+^*(\vec s)f \| =0,\label{lim: MTf}
\end{equation}
where the operator $\tilde\Omega_+^*(\vec s)$ is given by 
\begin{equation}
    \tilde \Omega_+^*(\vec s)=\Omega_+^*(t)
\end{equation}
for 
\begin{equation}
    s_j=t\, \quad\text{mod}\,\, T_j, \qquad j=1,\cdots,N.\label{eq: sandt}
\end{equation}
We note that for each $\vec s\in \bar{\mathbb T}_1\times \cdots\times \bar{\mathbb T}_N$, 
\begin{equation}
 \lim\limits_{M\to \infty}   \| (\mathbbm1-F_M(x,p))\tilde\Omega_+^*(\vec s)f \|=0.\label{lim: 1-FMf}
\end{equation}
By Lemma~\ref{lem: Omegas-uF} and~\eqref{lim: 1-FMf}, we obtain 
\begin{equation}
   \lim\limits_{M\to \infty} \lim\limits_{\vec u \to \vec s} \| F_M(x,p) \tilde \Omega_+^*(\vec u)\tilde P_c(\vec s)f\|= \lim\limits_{M\to \infty} \|F_M(x,p) \tilde \Omega_+^*(\vec s)f\|=\| \tilde P_c(\vec s)f\|,
\end{equation}
which yields 
\begin{equation}
   \lim\limits_{M\to \infty} \lim\limits_{\vec u \to \vec s} \| (\mathbbm1-F_M(x,p)) \tilde \Omega_+^*(\vec u)\tilde P_c(\vec s)f\|= 0.
\end{equation}
This together with~\eqref{lim: 1-FMf} yields 
\begin{equation}
\lim\limits_{M\to \infty}  \lim\limits_{\vec u\to \vec s}  \| (\mathbbm 1-F_M(x,p)) (\tilde \Omega_+^*(\vec u)-\tilde \Omega_+^*(\vec s))\tilde P_c(\vec s)f \|=0,
\end{equation}
where we used $\tilde \Omega_+^*(\vec s)=\tilde \Omega_+^*(\vec s)\tilde P_c(\vec s)$. This together with Proposition~\ref{prop: 4.2} yields 
\begin{equation}
   \lim\limits_{M\to \infty}  \lim\limits_{\vec u\to \vec s}  \| (\mathbbm 1-F_M(x,p)) \tilde \Omega_+^*(\vec u)f-(\mathbbm 1-F_M(x,p))\tilde \Omega_+^*(\vec s)f \|=0.
\end{equation}
Then this together with~\eqref{lim: 1-FMf}, compactness of $\bar{\mathbb T}_1\times \cdots \times \bar{\mathbb T}_N$ and the standard compactness argument, yields~\eqref{lim: MTf}. \end{proof}

    \subsection{Properties of the operators $\mathcal C_r(t)$}
In this section, we assume that $M\geq M_0$, where $M_0$ is given in Proposition~\ref{prop: 4.1}, and prove 
\begin{proposition}\label{BP1} If Assumption~\ref{aspV2} holds true, then 
\begin{equation}
   \left( \int_0^\infty \| \langle x\rangle^{-\eta} (\mathbbm1-\mathcal C_r(t))^{-1}P^\pm e^{-itH_0}f\|^2dt\right)^{1/2}\lesssim_{n,\eta}\|f\|, \qquad f\in \mathcal L^2_x(\mathbb R^n)
\end{equation}
holds true for all $\eta>\frac{3}{2}$ and $n\geq 8$.
\end{proposition}
and
\begin{proposition}\label{BP2} Under Assumptions~\ref{aspV2} and~\ref{aspV3}, we have
\begin{equation}
    \int_0^\infty\| \langle x\rangle^{-\eta}(\mathbbm 1-\mathcal C_r(t))^{-1}P^\pm e^{-itH_0}f \|dt\lesssim_{n,\eta}\|f\|\qquad f\in \mathcal L^2_x(\mathbb R^n)\label{est: BP2}
\end{equation}
for all $\eta>\frac{3}{2}$ and $n\geq 5$.
    
\end{proposition}
Let 
\begin{equation}
    I^\pm (t):=\int_0^\infty P^\pm e^{\pm is H_0}V(x,t+s)U(t+s,t) ds, \qquad t\in \mathbb R.\label{def: Ipmt}
\end{equation}
The proof of Propositions~\ref{BP1} and~\ref{BP2} uses the lemma and two propositions listed below. 
\begin{lemma}\label{lem: CM est} For all $M\geq 1$, we have
\begin{equation}
   \sup\limits_{t\in \mathbb R} \|\mathcal C_{M}(t) |p|^{-1/2}\langle x\rangle^{1/2+\epsilon}\|\lesssim_{M,\epsilon} 1\label{BP2: est1}
\end{equation}
and
\begin{equation}
   \sup\limits_{t\in \mathbb R} \|\mathcal C_{M}(t) \frac{\langle p\rangle^\alpha}{|p|^\alpha}\|\lesssim_{M,\alpha} 1\label{BP2: est2}
\end{equation}
for all $\epsilon\in (0,1/2)$ and $\alpha\in [0,\frac{7}{4})$ when $n\geq 4$. 
    
\end{lemma}
\begin{proof} \eqref{BP2: est1} follows by a proof similar to that of Proposition~\ref{prop: CM est}. To obtain~\eqref{BP2: est2}, we use the Duhamel's principle to obtain 
    \begin{align}
        \mathcal C_M(t)\frac{\langle p\rangle^\alpha}{|p|^\alpha}=&\mathcal C(t)\Omega_+^*(t) F_M(x,p)\Omega_+^*(t)\frac{\langle p\rangle^\alpha}{|p|^\alpha}=\mathcal C_{M,1}(t)+\mathcal C_{M,2}(t),\label{eq: decom CM1+CM2}
    \end{align}
where the operators $\mathcal C_{M,1}(t)$ and $\mathcal C_{M,2}(t)$ are given by 
\begin{equation}
        \mathcal C_{M,1}(t):= \mathcal C(t)\Omega_+^*(t) F_M(x,p)\frac{\langle p\rangle^\alpha}{|p|^\alpha}
\end{equation}
and
\begin{equation}
    \mathcal C_{M,2}(t):=(-i)\int_0^\infty \mathcal C(t)\Omega_+^*(t) F_M(x,p)e^{isH_0}V(x,t+s)U(t+s,t)\frac{\langle p\rangle^\alpha}{|p|^\alpha} ds.\label{def: CM2t2}
\end{equation}
By Eq.~\eqref{def: FMxp} and Corollary~\ref{cor: bd: C}, $\mathcal C_{M,1}(t)$ satisfies 
\begin{equation}
   \sup\limits_{t\in \mathbb R} \|\mathcal C_{M,1}(t)\|\lesssim_{M,\alpha} 1.\label{est: CM1t}
\end{equation}
To estimate $U(t+s,t)\frac{\langle p\rangle^\alpha}{|p|^\alpha}$, we note that by the Duhamel's principle, Hardy-Littlewood Sobolev inequality and the unitarity of $U(t+s,t+u)$, $(U(t+s,t)-e^{-isH_0})\frac{\langle p\rangle^\alpha}{|p|^\alpha}$ satisfies for $s\geq 0$,
\begin{align}
  \| (U(t+s,t)-e^{-isH_0})\frac{\langle p\rangle^\alpha}{|p|^\alpha}\|\leq &\int_0^s\| U(t+s,t+u)V(x,t+u)e^{-iuH_0}\frac{\langle p\rangle^\alpha}{|p|^\alpha} \|du\nonumber\\
  \lesssim & s\sup\limits_{u\in \mathbb R} \|\langle x\rangle^\alpha V(x,u) \|_{\mathcal L^\infty_x(\mathbb R^n)}.
\end{align}
This together with Eq.~\eqref{def: CM2t2}, Hardy-Littlewood Sobolev inequality and estimate 
\begin{equation}
    \|F_{M}(x,p)e^{isH_0}\langle x\rangle^{-\beta}\|_{2\to 2}\lesssim_{M, \beta}\frac{1}{\langle s\rangle^\beta},\qquad \beta>0, 
\end{equation}
yields 
\begin{align}
   & \|\mathcal C_{M,2}(t)\frac{\langle p\rangle^\alpha}{|p|^\alpha} \|\nonumber\\
   \leq & \int_0^\infty\| \mathcal C(t)\Omega_+^*(t)\|\|F_M(x,p)e^{isH_0}\langle x\rangle^{-5/4} \| \|\langle x\rangle^{3} V(x,t+s)\|_{\mathcal L^\infty_x(\mathbb R^n)}\|\langle x\rangle^{-\frac{7}{4}}\frac{\langle p\rangle^\alpha}{|p|^\alpha}e^{-isH_0}\|  ds\nonumber\\
     +&\int_0^\infty \| \mathcal C(t)\Omega_+^*(t)\|\|F_M(x,p)e^{isH_0}\langle x\rangle^{-9/4} \| \|\langle x\rangle^{9/4} V(x,t+s)\|_{\mathcal L^\infty_x(\mathbb R^n)}\| (U(t+s,t)-e^{-isH_0})\frac{\langle p\rangle^\alpha}{|p|^\alpha}\|  ds\nonumber\\
     \lesssim & \int_0^\infty  \frac{1}{\langle s\rangle^{5/4}} \sup\limits_{u\in \mathbb R}\|\langle x\rangle^{3} V(x,u)\|_{\mathcal L^\infty_x(\mathbb R^n)} ds+\int_0^\infty \frac{s}{\langle s\rangle^{9/4}} \|\langle x\rangle^{9/4} V(x,u) \|_{\mathcal L^\infty_x(\mathbb R^n)} ds\nonumber\\
     \lesssim & \|\langle x\rangle^3V(x,u) \|_{\mathcal L^\infty_x(\mathbb R^n)}
\end{align}
for $\alpha\in [0, \frac{7}{4})$. This together with~\eqref{est: CM1t} and Eq.~\eqref{eq: decom CM1+CM2} yields~\eqref{BP2: est2}.\end{proof}
\begin{proposition}\label{prop: main} Let $n\geq 5$ and $\alpha>3/2$ be two fixed numbers. If $\langle x\rangle^{\frac{9}{2}+\alpha} V(x,t)\in L^\infty_{x,t}(\mathbb R^{n+1})$, then for all $\epsilon\in (0, \min\{\alpha/2-3/4, 1/2\})$,
\begin{equation}
    \| \frac{|p|^\alpha}{\langle p\rangle^\alpha} I^\pm(t) |p|^{-1/2}\langle x\rangle^{1/2+\epsilon}\|\leq C\max\limits_{j=1,2} \| \langle x\rangle^{\frac{9}{2}+\alpha} V(x,s)\|_{L^\infty_{x,s}(\mathbb R^{n+1})}^j,\label{B: prop: est}
\end{equation}
where the constant $C=C(\epsilon,n,\alpha)>0$ stands for a constant depending on $\alpha, n$ and $\epsilon$.
\end{proposition}
\begin{proposition}\label{prop: main2} There exists $\beta\in (0, 1/8]$ such that for all $\alpha\in (0, \frac{7}{4})$, with $\alpha\beta_+:=\max\{\alpha-\beta, 0\}$, we have
\begin{equation}
    \| \frac{|p|^{\alpha\beta_+}}{\langle p\rangle^{\alpha\beta_+}}I^\pm (t)\frac{\langle p\rangle^{\alpha}}{|p|^{\alpha}}\|\leq C ,
\end{equation}
where the constant $C=C(\epsilon,n)$ stands for a constant depending on $n$ and $\epsilon$.

\end{proposition}
We prove Proposition~\ref{prop: main} first. The proof of Proposition~\ref{prop: main} requires the following lemma.
\begin{lemma}\label{B:lem: main}If $\langle x\rangle^{\sigma}V(x,t)\in L^\infty_{x,t}(\mathbb R^{n+1})$ for some $\sigma\geq 2$, then for all $n\geq 5$ and $\epsilon\in (0,1/3)$,
\begin{equation}
    \int_0^\infty \langle s\rangle^{-2-3\epsilon}\| \langle x\rangle^{-2}U(t+s,t)|p|^{-1/2}\langle x\rangle^{1/2+\epsilon}\|ds\leq C\max\limits_{j=0,1} \| \langle x\rangle^\sigma V(x,s)\|_{L^\infty_{x,s}(\mathbb R^{n+1})}^j,\label{BLem5: est}
\end{equation}
where $C=C(n,\epsilon)>0$ stands for a constant depending on $n$ and $\epsilon$.
\end{lemma}
We defer the proof of Lemma~\ref{B:lem: main} after the proof of Proposition~\ref{prop: main}.

\begin{proof}[Proof of Proposition~\ref{prop: main}] We estimate $ \| \frac{|p|^\alpha}{\langle p\rangle^\alpha} I^+(t) |p|^{-1/2}\langle x\rangle^{1/2+\epsilon}\|$ and $ \| \frac{|p|^\alpha}{\langle p\rangle^\alpha} I^-(t) |p|^{-1/2}\langle x\rangle^{1/2+\epsilon}\|$ can be treated similarly. By estimate~\eqref{Sep20.1}, we obtain, for all $\epsilon_1>0$ satisfying condition $\frac{5}{4}+\alpha/2-\epsilon_1>2$, 
\begin{align}
   & \| \frac{|p|^\alpha}{\langle p\rangle^\alpha} I^+(t) |p|^{-1/2}\langle x\rangle^{1/2+\epsilon}\|\nonumber\\
   \leq & \int_0^\infty  \|\frac{|p|^\alpha}{\langle p\rangle^\alpha} P^+ e^{isH_0} \langle x\rangle^{-(\frac{5}{2}+\alpha)}\| \| \langle x\rangle^{\frac{9}{2}+\alpha} V(x,u)\|_{L^\infty_{x,u}(\mathbb R^{n+1})} \| \langle x\rangle^{-2}U(t+s,t) |p|^{-1/2}\langle x\rangle^{1/2+\epsilon}\| ds\nonumber\\
   \lesssim_{\epsilon_1,\alpha}  &  \| \langle x\rangle^{\frac{9}{2}+\alpha} V(x,u)\|_{L^\infty_{x,u}(\mathbb R^{n+1})}\int_0^\infty \frac{1}{\langle s\rangle^{\frac{5}{4}+\alpha/2-\epsilon_1}} \| \langle x\rangle^{-2}U(t+s,t) |p|^{-1/2}\langle x\rangle^{1/2+\epsilon}\| ds.
\end{align}
This together with Lemma~\ref{B:lem: main} (take $\epsilon=\frac{1}{3}(-\frac{3}{4}+\frac{\alpha}{2}-\epsilon_1)>0$) implies~\eqref{B: prop: est}. \end{proof}
Now we prove Lemma~\ref{B:lem: main}. The proof of Lemma~\ref{B:lem: main} requires 
\begin{lemma}\label{B: lem1}For all $n\geq 5$,
    \begin{equation}
    \| \langle x\rangle^{-2}|p|^{-1/2}\langle x\rangle^{3/2}\|\lesssim_n 1.\label{B7}
\end{equation}
\end{lemma}
\begin{proof} Let 
\begin{equation}
   E:= \| \langle x\rangle^{-2} |p|^{-1/2} \langle x\rangle^{3/2}\|.
\end{equation}
We write $\langle x\rangle^{3/2}$ as $\langle x\rangle^2\langle x\rangle^{-1/2}$ and note that
\begin{equation}
  E \leq \| \langle x\rangle^{-2} [|p|^{-1/2},  \langle x\rangle^2] \langle x\rangle^{-1/2} \|+\|  |p|^{-1/2}\langle x\rangle^{-1/2}\|.\label{B9}
\end{equation}
By Hardy-Littlewood-Sobolev inequality, 
\begin{align}
    \| |p|^{-1/2}\langle x\rangle^{-1/2}\|\leq \||p|^{-1/2}|x|^{-1/2}\||x|^{1/2}\langle x\rangle^{-1/2}\|\lesssim 1,\label{B10}
\end{align}
for all $n\geq 2$. To estimate $\| \langle x\rangle^{-2} [|p|^{-1/2},  \langle x\rangle^2] \langle x\rangle^{-1/2} \|$, we also note that 
\begin{align}
 &\| \langle x\rangle^{-2}   [|p|^{-1/2},  \langle x\rangle^2]\langle x\rangle^{-1/2}\|\nonumber\\
 \leq& \sum\limits_{j=1}^n \left(\| \langle x\rangle^{-2} \partial^2_{p_j}[|p|^{-1/2}] \langle x\rangle^{-1/2}\|+2 \| x_j\langle x\rangle^{-2} \partial_{p_j}[|p|^{-1/2}] \langle x\rangle^{-1/2}\|\right),
\end{align}
which implies 
\begin{align}
  &  \| \langle x\rangle^{-2}   [|p|^{-1/2},  \langle x\rangle^2]\langle x\rangle^{-1/2}\|\nonumber\\
  \lesssim &\sum\limits_{j=1}^n \left( \|\langle x\rangle^{-2} |p|^{-5/2}\langle x\rangle^{-1/2} \| +\|\langle x\rangle^{-2} p_j^2|p|^{-9/2}\langle x\rangle^{-1/2} \| \right)\nonumber\\
  &+\sum\limits_{j=1}^n \|x_j\langle x\rangle^{-2} p_j|p|^{-5/2} \langle x\rangle^{-1/2} \|.
\end{align}
By employing estimate $\|p_j/|p|\|\leq 1$ and Hardy-Littlewood-Sobolev inequality, we have, for $n\geq 5$,
\begin{align}
    \|\langle x\rangle^{-2} |p|^{-5/2}\langle x\rangle^{-1/2} \|\leq& \| |x|^2\langle x\rangle^{-\sigma/2}\|\| |x|^{-2}|p|^{-2}\|\||p|^{-1/2}|x|^{-1/2}\|\||x|^{1/2}\langle x\rangle^{-1/2}\|\nonumber\\
    \lesssim &1,
\end{align}
\begin{align}
    \|\langle x\rangle^{-2} p_j^2|p|^{-9/2}\langle x\rangle^{-1/2} \|\leq & \| |x|^2\langle x\rangle^{-2}\|\| |x|^{-2}|p|^{-2}\| \|p_j/|p|\|^2\||p|^{-1/2}|x|^{-1/2}\|\||x|^{1/2}\langle x\rangle^{-1/2}\|\nonumber\\
    \lesssim & 1
\end{align}
and
\begin{align}
    \|x_j\langle x\rangle^{-2} p_j|p|^{-5/2} \langle x\rangle^{-1/2} \|\leq& \|x_j\langle x\rangle^{-2}|x| \|\| |x|^{-1}|p|^{-1}\|\| p_j/|p|\|\||p|^{-1/2}|x|^{-1/2}\|\||x|^{1/2}\langle x\rangle^{-1/2}\|\nonumber\\
    \lesssim &1.
\end{align}
These yield 
\begin{equation}
    \|\langle x\rangle^{-2}[|p|^{-1/2}, \langle x\rangle^2]\langle x\rangle^{-1/2}\|\lesssim_n 1.
\end{equation}
This together with~\eqref{B9} and~\eqref{B10} yields~\eqref{B7}.\end{proof}
\begin{lemma}\label{B:lem2} For all $n\geq 4$, 
    \begin{equation}
  \|\langle x\rangle^{-3/2} |p|^{-1/2} x_j\|\lesssim 1, \qquad j=1,\cdots,n.\label{Blem: est}
\end{equation}

\end{lemma}
\begin{proof} We note that 
\begin{align}
     \|\langle x\rangle^{-3/2} |p|^{-1/2} x_j\|\leq & \|x_j \langle x\rangle^{-3/2} |p|^{-1/2}\|+\|\langle x\rangle^{-3/2}[x_j, |p|^{-1/2}]\|\nonumber\\
     =& \|x_j \langle x\rangle^{-3/2} |p|^{-1/2}\|+\frac{1}{2}\| \langle x\rangle^{-3/2} p_j|p|^{-5/2}\|.\label{Blem: eq1}
\end{align}
By Hardy-Littlewood-Sobolev inequality and $\|p_j/|p|\|\leq 1$, we obtain 
\begin{equation}
    \|x_j \langle x\rangle^{-3/2} |p|^{-1/2}\|\leq \|x_j \langle x\rangle^{-3/2}|x|^{1/2} \|\| |x|^{-1/2}|p|^{-1/2}\|\lesssim 1
\end{equation}
and
\begin{equation}
    \|\langle x\rangle^{-3/2} p_j|p|^{-5/2}\|\leq \| \langle x\rangle^{-3/2}|x|^{3/2} \|\| |x|^{-3/2}|p|^{-3/2}\|\|p_j/|p|\|\lesssim 1.
\end{equation}
These together with~\eqref{Blem: eq1} imply~\eqref{Blem: est}.\end{proof}
\begin{lemma}\label{B:lem4} For all $n\geq 3$ and $\epsilon\in (0,1/2)$,
    \begin{equation}
 \int   \langle s\rangle^{-1-2\epsilon} \| \langle x\rangle^{-2} e^{-isH_0} |p|^{-1/2} \langle x\rangle^{1/2+\epsilon} \| ds\leq C,\label{B21}
\end{equation}
where $C=C(n,\epsilon)>0$ stands for a constant depending on $n$ and $\epsilon$.
\end{lemma}
\begin{proof} By endpoint Strichartz estimate, Lemma~\ref{B: lem1} and H\"older's inequality, we have
\begin{align}
   & \int   \langle s\rangle^{-1-2\epsilon} \| \langle x\rangle^{-2} e^{-isH_0} |p|^{-1/2} \langle x\rangle^{1/2+\epsilon} \chi(|x|\leq \langle s\rangle)\| ds\nonumber\\
   \lesssim & \int   \langle s\rangle^{-1-2\epsilon} \| \langle x\rangle^{-2} |p|^{-1/2} \langle x\rangle^{3/2} \|\|\langle x\rangle^{-3/2} e^{-isH_0} \| \langle s\rangle^{1/2+\epsilon} ds\nonumber\\
   \lesssim & \| \langle s\rangle^{-1/2-\epsilon}\|_{L^2_s(\mathbb R)} \| e^{-isH_0}\|_{L^2_x(\mathbb R^n)\to L^2_sL^{\frac{2n}{n-2}}_x(\mathbb R^{n+1})}\nonumber\\
   \lesssim_{\epsilon,n} & 1\label{BQs}
\end{align}
for $n\geq 3$ and $\epsilon\in (0,1/2)$. Next, we estimate 
\begin{equation}
    Q:= \int   \langle s\rangle^{-1-2\epsilon} \| \langle x\rangle^{-2} e^{-isH_0} |p|^{-1/2} \langle x\rangle^{1/2+\epsilon}\chi(|x|>\langle s\rangle) \| ds.
\end{equation}
To estimate~$Q$, we let $\{F_j\}_{j=1}^{j=n}$ denote a partition of the unity satisfying 
\begin{equation}
    |x_j|\geq \frac{|x|}{n}, \qquad \text{for all }x_j \in \text{supp}(F_j).
\end{equation}
This implies 
\begin{equation}
    F_j\chi(|x|>\langle s\rangle)=F_j\chi(|x|>\langle s\rangle)\chi(|x_j|>\frac{\langle s\rangle}{n})
\end{equation}
and
\begin{equation}
     \| F_j\langle x_j\rangle^{-1/2-\epsilon}\langle x\rangle^{1/2+\epsilon}\chi(|x|>1)\|\lesssim_n 1 \label{Fj, xp}.
\end{equation}
By triangle inequality and estimate~\eqref{Fj, xp}, we find 
\begin{equation}
    Q\lesssim_n \sum\limits_{j=1}^n Q_j,\label{BQnQj}
\end{equation}
where $Q_j, j=1,\cdots,n$, are defined by 
\begin{equation}
    Q_j:=\int  \langle s\rangle^{-1-2\epsilon} \| \langle x\rangle^{-2} e^{-isH_0} |p|^{-1/2} \langle x_j\rangle^{1/2+\epsilon}\chi(|x_j|>\frac{\langle s\rangle}{n}) \| ds.
\end{equation}
We write 
\begin{equation}
    e^{-isH_0} x_j=(x_j-2sp_j)e^{-isH_0}
\end{equation}
to obtain 
\begin{equation}
    Q_j\leq Q_{j1}+Q_{j2}+Q_{j3},\label{B29}
\end{equation}
where 
\begin{equation}
    Q_{j1}:=\int  \langle s\rangle^{-1-2\epsilon} \| \langle x\rangle^{-2}x_j|p|^{-1/2} e^{-isH_0}  \frac{\langle x_j\rangle^{1/2+\epsilon}}{x_j}\chi(|x_j|>\frac{\langle s\rangle}{n}) \| ds,
\end{equation}
\begin{equation}
    Q_{j2}:=2\int   |s|\langle s\rangle^{-1-2\epsilon} \| \langle x\rangle^{-2} e^{-isH_0} |p|^{-1/2}p_j \frac{\langle x_j\rangle^{1/2+\epsilon}}{x_j}\chi(|x_j|>\frac{\langle s\rangle}{n}) \| ds
\end{equation}
and
\begin{equation}
    Q_{j3}:=\int  \langle s\rangle^{-1-2\epsilon} \| \langle x\rangle^{-2} e^{-isH_0} [x_j,|p|^{-1/2}] \frac{\langle x_j\rangle^{1/2+\epsilon}}{x_j}\chi(|x_j|>\frac{\langle s\rangle}{n}) \| ds.
\end{equation}
By employing estimate
\begin{equation}
    \| \frac{\langle x_j\rangle^{1/2+\epsilon}}{x_j}\chi(|x_j|>\frac{\langle s\rangle}{n})\|\lesssim_n \frac{1}{\langle s\rangle^{1/2-\epsilon}},\qquad \forall\, \epsilon\in (0,1/2),\label{we: ns}
\end{equation}
the unitarity of $e^{-isH_0}$ and Hardy-Littlewood Sobolev inequality, we have for all $n\geq 3$,
\begin{align}
    Q_{j1}\leq & \int   \langle s\rangle^{-1-2\epsilon} \| \langle x\rangle^{-2}x_j|p|^{-1/2} \|\|e^{-isH_0}\|\|  \frac{\langle x_j\rangle^{1/2+\epsilon}}{x_j}\chi(|x_j|>\frac{\langle s\rangle}{n}) \| ds\nonumber\\
    \lesssim & \int \langle s\rangle^{-1-\epsilon} ds  \lesssim_\epsilon 1.\label{B: Qj1}
\end{align}
By employing H\"older's inequality and estimates~\eqref{we: ns}, $\| \frac{p_j}{|p|}\|\leq 1$ and 
\begin{equation}
    \int \| \langle x\rangle^{-2} |p|^{1/2}e^{-isH_0}f\|^2 ds \lesssim \|f\|^2,\qquad  f\in L^2(\mathbb R^n),
\end{equation}
we obtain 
\begin{align}
    Q_{j2}\leq & 2 \int   |s|\langle s\rangle^{-1-2\epsilon} \| \langle x\rangle^{-2} e^{-isH_0} |p|^{1/2}\|\|\frac{p_j}{|p|}\|\| \frac{\langle x_j\rangle^{1/2+\epsilon}}{x_j}\chi(|x_j|>\frac{\langle s\rangle}{n}) \| ds\nonumber\\
    \lesssim_n &  \int  \langle s\rangle^{-1/2-\epsilon} \|\langle x\rangle^{-2} e^{-isH_0} |p|^{1/2} \| ds \nonumber\\
    \lesssim_{n} & \| \langle s\rangle^{-1/2-\epsilon} \|_{L^2_s(\mathbb R)}\| \langle x\rangle^{-2} e^{-isH_0} |p|^{1/2} \|_{L^2_x(\mathbb R^n)\to L^2_{x,s}(\mathbb R^{n+1})}\nonumber\\
    \lesssim_{n,\epsilon} & 1.\label{B: Qj2}
\end{align}
Using $[x_j,|p|^{-1/2}]=\frac{cp_j}{|p|^{5/2}}$ for some constant $c>0$ and Hardy-Littlewood Sobolev inequality, we estimate
\begin{align}
    Q_{j3}\leq &c\int   \langle s\rangle^{-1-2\epsilon} \| \langle x\rangle^{-2} e^{-isH_0} \frac{p_j}{|p|^{5/2}}\frac{\langle x_j\rangle^{1/2+\epsilon}}{x_j}\chi(|x_j|>\frac{\langle s\rangle}{n}) \| ds\nonumber\\
    \leq & c\int   \langle s\rangle^{-1-2\epsilon} \|\langle x\rangle^{-2}|p|^{-1-\epsilon} \|\|e^{-isH_0}\|\|\frac{p_j}{|p|}\|\||p|^{-1/2+\epsilon}\langle x\rangle^{-1/2+\epsilon} \|ds\nonumber\\
    \lesssim_\epsilon & 1.
\end{align}
This together with estimates~\eqref{B29},~\eqref{B: Qj1} and~\eqref{B: Qj2} yields 
\begin{equation}
    Q_j\lesssim_{n,\epsilon} 1,\qquad j=1,\cdots,n.\label{B: est: Qj}
\end{equation}
Estimates~\eqref{BQs},~\eqref{BQnQj} and~\eqref{B: est: Qj} yield~\eqref{B21}. 
\end{proof}
\begin{proof}[Proof of Lemma~\ref{B:lem: main}]By Lemma~\ref{B:lem4}, it suffices to estimate 
\begin{equation}
    Q:=\int_0^\infty \langle s\rangle^{-2-3\epsilon}\| \langle x\rangle^{-2}(U(t+s,t)-e^{-isH_0})|p|^{-1/2}\langle x\rangle^{1/2+\epsilon}\|ds.\label{Blem:5Q}
\end{equation}
By the Duhamel principle, unitarity of $U(t+s,t+u)$ and Lemma~\ref{B:lem4}, we have 
\begin{align}
 & \|  (U(t+s,t)-e^{-isH_0})|p|^{-1/2}\langle x\rangle^{1/2+\epsilon}\|\nonumber\\
 \leq &\int_0^s \|U(t+s,t+u)V(x,t+u)e^{-iuH_0}|p|^{-1/2}\langle x\rangle^{1/2+\epsilon}\| du  \nonumber\\
 \leq & \| \langle x\rangle^2V(x,v)\|_{L^\infty_{x,v}(\mathbb R^{n+1})} \langle s\rangle^{1+2\epsilon}\int_0^s \langle u\rangle^{-1-2\epsilon}\| \langle x\rangle^{-2} e^{-iuH_0}|p|^{-1/2}\langle x\rangle^{1/2+\epsilon}\|du\nonumber\\
 \leq & C\| \langle x\rangle^2V(x,v)\|_{L^\infty_{x,v}(\mathbb R^{n+1})} \langle s\rangle^{1+2\epsilon},
\end{align}
where $C=C(n,\epsilon)>0$ stands for a constant depending on $n$ and $\epsilon$. This together with Eq.~\eqref{Blem:5Q} yields 
\begin{equation}
    Q\leq C\int_0^\infty  \langle s\rangle^{-1-\epsilon} \| \langle x\rangle^2V(x,v)\|_{L^\infty_{x,v}(\mathbb R^{n+1})}  ds\leq C_0\| \langle x\rangle^2V(x,v)\|_{L^\infty_{x,v}(\mathbb R^{n+1})} , 
\end{equation}
where $C_0=C_0(n,\epsilon)>0$ stands for a constant depending on $n$ and $\epsilon$. This together with Lemma~\ref{B:lem4} yields~\eqref{BLem5: est}.\end{proof}
Now we prove Proposition~\ref{prop: main2}. The proof of Proposition~\ref{prop: main2} requires
\begin{lemma}\label{lem: prop: main2}For all $\alpha\in (0, \frac{7}{4})$ and $n\geq5$,
\begin{equation}
  \| \langle x\rangle^{-\alpha}U(t+s,t) \frac{\langle p\rangle^{\alpha}}{|p|^\alpha}\|\lesssim_n \langle s\rangle^{\alpha/2}
\end{equation}
holds true.
\end{lemma}
\begin{remark} The upper bound $\frac{7}{4}$ for $\alpha$ is not sharp but sufficient for Proposition~\ref{prop: main2}.
    
\end{remark}
\begin{proof} Let $\mathcal B(t,s)\equiv\langle x\rangle^{-\alpha}U(t+s,t) \frac{\langle p\rangle^{\alpha}}{|p|^\alpha}$. We deal with the case when $s\geq 0$. The case when $s<0$ can be treated similarly. We break $\mathcal B(t,s)$ into two pieces
\begin{equation}
    \mathcal B(t,s)=\mathcal B_1(t,s)+\mathcal B_2(t,s),\label{id: Bts}
\end{equation}
where the operators $\mathcal B_j(t,s), j=1,2,$ are given by 
\begin{equation}
    \mathcal B_1(t,s):=\mathcal B(t,s)\chi(|p|> \frac{1}{\langle s\rangle^{1/2}}) 
\end{equation}
and $\mathcal B_2(t,s):=\mathcal B(t,s)-\mathcal B_1(t,s)$. Using 
\begin{equation}
  \| \frac{\langle p\rangle^{\alpha}}{|p|^\alpha}\chi(|p|> \frac{1}{\langle s\rangle^{1/2}}) \|\lesssim \langle s\rangle^{\alpha/2 }
\end{equation}
and the unitarity of $U(t+s,t)$, we obtain 
\begin{equation}
   \| \mathcal B_1(t,s) \frac{\langle p\rangle^\alpha}{|p|^\alpha} \|\lesssim \langle s\rangle^{\alpha/2}.\label{est: B1ts}
\end{equation}
Next, we estimate $\mathcal B_2(t,s)\frac{\langle p\rangle^\alpha}{|p|^{\alpha}}$. Using
\begin{equation}
    \| \langle x\rangle^{-\alpha} \frac{\langle p\rangle^\alpha}{|p|^\alpha} e^{-iuH_0}\chi(|p|\leq \frac{1}{\langle s\rangle^{1/2}}) \|\lesssim 1,\qquad \forall u\in \mathbb R,\, \alpha\in [0,\frac{7}{4}),
\end{equation}
by the Duhamel's principle, the unitarity of $U(t+s,t+u)$ and Hardy-Littlewood-Sobolev inequality, we arrive at, for $s\geq 0$, 
\begin{align}
    \| \mathcal B_2(t,s)\frac{\langle p\rangle^\alpha}{|p|^\alpha}\|\leq &  \| \langle x\rangle^{-\alpha}e^{-isH_0}\frac{\langle p\rangle^\alpha}{|p|^\alpha}\chi(|p|\leq \frac{1}{\langle s\rangle^{1/2}})\|\nonumber\\
    &+\int_0^s \|\langle x\rangle^{-\alpha} U(t+s,t+u)V(x,t+u)e^{-iuH_0} \frac{\langle p\rangle^\alpha}{|p|^\alpha}\chi(|p|\leq \frac{1}{\langle s\rangle^{1/2}}) \|du\nonumber\\
    \lesssim &  1+\int_0^s\sup\limits_{u\in \mathbb R} \| \langle x\rangle^{2} V(x,u)\|_{\mathcal L^\infty_x(\mathbb R^n)}\||x|^{-2}|p|^{-2}\|\| \langle p\rangle^\alpha |p|^{2-\alpha}\chi(|p|\leq \frac{1}{\langle s\rangle^{1/2}})\|du\nonumber\\
    \lesssim & \langle s\rangle^{\alpha/2}.
\end{align}
This together with estimate~\eqref{est: B1ts} and Eq.~\eqref{id: Bts} yields 
\begin{equation}
    \|\mathcal B(t,s)\frac{\langle p\rangle^\alpha}{|p|^\alpha}\|\lesssim \langle s\rangle^{\alpha/2}.
\end{equation}

\end{proof}
Now we prove Proposition~\ref{prop: main2}.
\begin{proof}[Proof of Proposition~\ref{prop: main2}] By Lemma~\ref{lem: prop: main2} and estimate~\eqref{Sep20.1}, we obtain for $\alpha\in (0, \frac{7}{4})$,
\begin{equation}
    \| \frac{|p|^{\alpha\beta_+}}{\langle p\rangle^{\alpha\beta_+}}I^\pm(t)\frac{\langle p\rangle^\alpha}{|p|^{\alpha}}\|\lesssim_{\epsilon, n} \int_0^\infty \frac{1}{\langle s\rangle^{5/4+\alpha\beta_+/2-\epsilon}}\sup\limits_{u\in \mathbb R}\|\langle x\rangle^{2\alpha+5/2}V(x,u)\|_{\mathcal L^\infty_x(\mathbb R^n)} \langle s\rangle^{\alpha/2}ds.\label{pr: lem: est1}
\end{equation}
Then with $\beta=\frac{1}{8}$ and $\epsilon\in (0, 1/8)$, we obtain 
\begin{align}
    5/4+\alpha\beta_+/2-\epsilon -\alpha/2 =&5/4+\max\{\alpha-\beta, 0\}/2-\epsilon-\alpha/2\nonumber\\
    >& 5/4-\frac{1}{16}-\frac{1}{8}>1. 
\end{align}
This together with~\eqref{pr: lem: est1} yields for all $\epsilon\in (0, \frac{1}{8})$,
\begin{equation}
    \| \frac{|p|^{\alpha\beta_+}}{\langle p\rangle^{\alpha\beta_+}}I^\pm(t)\frac{\langle p\rangle^\alpha}{|p|^{\alpha}}\|\lesssim_{\epsilon, n} 1.
\end{equation}

\end{proof}
The proof of Proposition~\ref{BP1} also requires the following proposition.
\begin{proposition}\label{prop: weight}If Assumption~\ref{aspV2} is satisfied, then for all $\delta \in [0, \min\{\frac{n}{2}-2, 2\})$, $n\geq 5$,
    \begin{equation}
      \sup\limits_{t\in \mathbb R}  \| P_b(t)\langle x\rangle^\delta\|\lesssim_{\delta,n} 1.\label{est: Pbw}
    \end{equation}
\end{proposition}
\begin{proof}It suffices to fix $\delta \in [\max\{0,\frac{n}{2}+4-\sigma\}, \min\{\frac{n}{2}-2, 2\})$. It suffices to show that 
\begin{equation}
      \sup\limits_{t\in \mathbb R}  \| P_b(t)\frac{x_j^2}{\langle x\rangle^{2-\delta}}\|\lesssim_{\delta,n} 1,\qquad j=1,\cdots,n.
    \end{equation}
We estimate the case when $j=1$ and the case when $j=2,\cdots,n$, can be treated similarly. To estimate $ P_b(t)\frac{x_1^2}{\langle x\rangle^{2-\delta}}$, we break $P_b(t)\frac{x_1^2}{\langle x\rangle^{2-\delta}}$ into two parts 
\begin{equation}
    P_b(t)\frac{x_1^2}{\langle x\rangle^{2-\delta}}=P_b(t)x_1^2P^+\frac{1}{\langle x\rangle^{2-\delta}}+P_b(t)x_1^2P^-\frac{1}{\langle x\rangle^{2-\delta}}.\label{def: P++P-Pb}
\end{equation}
We note that with $P_b(t)\Omega_+(t)=0$, by the Duhamel's principle, the operator $P_b(t)x_1^2P^+\frac{1}{\langle x\rangle^{2-\delta}}$ reads 
\begin{align}
    P_b(t)x_1^2P^+ \frac{1}{\langle x\rangle^{2-\delta}}=&  P_b(t)(\mathbbm1 -\Omega_+(t))x_1^2P^+ \frac{1}{\langle x\rangle^{2-\delta}}\nonumber\\
    =&(-i)\int_0^\infty P_b(t)U(t,t+s)V(x,t+s)e^{-isH_0}x_1^2P^+\frac{1}{\langle x\rangle^{2-\delta}} ds.
\end{align}
Moving $x_1^2$ through $e^{-isH_0}$, we obtain 
\begin{align}
     P_b(t)x_1^2P^+ \frac{1}{\langle x\rangle^{2-\delta}}=&(-i)\int_0^\infty P_b(t)U(t,t+s)V(x,t+s)x_1^2e^{-isH_0}P^+\frac{1}{\langle x\rangle^{2-\delta}} ds\nonumber\\
     &+4i\int_0^\infty sP_b(t)U(t,t+s)V(x,t+s)x_1p_1e^{-isH_0}P^+\frac{1}{\langle x\rangle^{2-\delta}} ds\nonumber\\
     &+4(-i)\int_0^\infty s^2P_b(t)U(t,t+s)V(x,t+s)p_1^2e^{-isH_0}P^+\frac{1}{\langle x\rangle^{2-\delta}} ds.
\end{align}
By estimates~$\|P_b(t)\|\leq 1$,~\eqref{Oct.1},~\eqref{Sep20.1} and~\eqref{est: Pbweight} and the unitarity of $U(t,t+s)$, this implies 
\begin{align}
    &\| P_b(t)x_1^2P^+ \frac{1}{\langle x\rangle^{2-\delta}}\|\nonumber\\
    \leq & \int_0^\infty \|\langle x\rangle^5 V(x,t+s)\|_{\mathcal L^\infty_x(\mathbb R^n)}\|\langle x\rangle^{-3}e^{-isH_0}P^+\| ds\nonumber\\
    &+4\sum\limits_{l=1}^2\int_1^\infty s^l \|\langle x\rangle^{\frac{n}{2}+l+2-\delta} V(x,t+s)\|_{\mathcal L^\infty_x(\mathbb R^n)}\| \langle x\rangle^{-(\frac{n}{2}+l+2-\delta)}x_1^{2-l}p_1^le^{-isH_0}P^+\frac{1}{\langle x\rangle^{2-\delta}}\|ds\nonumber\\
     &+4\sum\limits_{l=1}^2\int_0^1 s^l \|\langle x\rangle^{\frac{n}{2}+l+2-\delta} V(x,t+s)\|_{\mathcal L^\infty_x(\mathbb R^n)}\| \langle x\rangle^{-(\frac{n}{2}+l+2-\delta)}x_1^{2-l}p_1^le^{-isH_0}P^+\frac{1}{\langle x\rangle^{2-\delta}}\|ds\nonumber\\
    \lesssim_{\epsilon,n} & \sup\limits_{u\in \mathbb R}\|\langle x\rangle^{\frac{n}{2}+4-\delta} V(x,u)\|_{\mathcal L^\infty_x(\mathbb R^n)}\left(\int_1^\infty \left(\frac{1}{\langle s\rangle^{3/2-\epsilon}}+\frac{1}{\langle s\rangle^{\frac{n}{4}+\frac{2-l-\delta}{2}-\epsilon}}\right) ds+1\right)\nonumber\\
    \lesssim_{\epsilon,n} & \sup\limits_{u\in \mathbb R}\|\langle x\rangle^{\frac{n}{2}+4-\delta} V(x,u)\|_{\mathcal L^\infty_x(\mathbb R^n)}
\end{align}
for all $\epsilon\in (0,\frac{n}{4}-\frac{\delta}{2}-1)$. This yields with $\epsilon=\frac{n}{8}-\frac{\delta}{4}-\frac{1}{2}$, 
\begin{equation}
     \sup\limits_{t\in \mathbb R} \| P_b(t)x_1^2P^+ \frac{1}{\langle x\rangle^{2-\delta}}\|\lesssim_{\delta,n} \sup\limits_{u\in \mathbb R}\|\langle x\rangle^{\frac{n}{2}+4-\delta} V(x,u)\|_{\mathcal L^\infty_x(\mathbb R^n)}.\label{est: P+t}
\end{equation}
By~\cite{JL1991}, we note that 
\begin{equation}
    P_b(t)=s\text{-}\lim\limits_{s\to -\infty} U(t,t+s) F_c(\frac{|x-2sp|}{s^\alpha}\geq 1)U(t+s,t),\qquad \text{ on }\mathcal L^2_x(\mathbb R^n)
\end{equation}
for all $\alpha\in (0, 1-\frac{2}{n}), \, n\geq 3.$ Hence, similarly, we have 
\begin{equation}
     \sup\limits_{t\in \mathbb R} \| P_b(t)x_1^2P^- \frac{1}{\langle x\rangle^{2-\delta}}\|\lesssim_{\delta,n} \sup\limits_{u\in \mathbb R}\|\langle x\rangle^{\frac{n}{2}+4-\delta} V(x,u)\|_{\mathcal L^\infty_x(\mathbb R^n)}.
\end{equation}
This together with~\eqref{est: P+t} and Eq.~\eqref{def: P++P-Pb} yields 
\begin{equation}
    \sup\limits_{t\in \mathbb R} \|P_b(t)\frac{x_1^2}{\langle x\rangle^{2-\delta}} \|\lesssim_{\delta,n} \sup\limits_{u\in \mathbb R}\|\langle x\rangle^{\frac{n}{2}+4-\delta} V(x,u)\|_{\mathcal L^\infty_x(\mathbb R^n)}.
\end{equation}
Similarly, we obtain 
\begin{equation}
    \sup\limits_{t\in \mathbb R} \|P_b(t)\frac{x_j^2}{\langle x\rangle^{2-\delta}} \|\lesssim_{\delta,n} \sup\limits_{u\in \mathbb R}\|\langle x\rangle^{\frac{n}{2}+4-\delta} V(x,u)\|_{\mathcal L^\infty_x(\mathbb R^n)}, \qquad j=2,\cdots,n.
\end{equation}
Thus, we conclude~\eqref{est: Pbw}.\end{proof}
Next, we prove Proposition~\ref{BP1}. 
\begin{proof}[Proof of Proposition~\ref{BP1}] Take $f\in \mathcal L^2_x(\mathbb R^n)$. We estimate $\langle x\rangle^{-\eta}(\mathbbm1-\mathcal C_r(t))^{-1}P^+ e^{-itH_0}f$ and $\langle x\rangle^{-\eta}(\mathbbm1-\mathcal C_r(t))^{-1}P^- e^{-itH_0}f$ can be treated similarly. We write $(\mathbbm1-\mathcal C_r(t))^{-1}P^+ e^{-itH_0}f $ as 
\begin{equation}
    (\mathbbm1-\mathcal C_r(t))^{-1}P^+ e^{-itH_0}f=f_1(t)+f_2(t),
\end{equation}
where $f_j(t), j=1,2,$ are given by 
\begin{equation}
    f_1(t):=P^+e^{-itH_0}f
\end{equation}
and 
\begin{equation}
    f_2(t):=(\mathbbm1-\mathcal C_r(t))^{-1}\mathcal C_r(t)P^+ e^{-itH_0}f.
\end{equation}
By estimate~\eqref{Nov.20.6} and the $L^2$ local decay of the free flow, $f_1(t)$ satisfies 
\begin{align}
    \left(\int_0^\infty \| \langle x\rangle^{-\eta}f_1(t)\|^2 dt\right)^{1/2}\leq & \left(\int_0^\infty \| \langle x\rangle^{-\eta}P^+\langle x\rangle^{3/2}\|^2\| \langle x\rangle^{-3/2}e^{-itH_0}f\|^2 dt\right)^{1/2}\lesssim \|f\|
\end{align}
for all $\eta>\frac{3}{2}$. For $f_2(t)$, we write $\mathcal C_r(t)=\mathcal C(t)-\mathcal C(t)P_b(t)-\mathcal C_M(t)$ to decompose $f_2(t)$ further 
\begin{equation}
    f_2(t)=\sum\limits_{j=1}^3 f_{2j}(t),
\end{equation}
where $f_{2j}(t), j=1,2,3,$ read 
\begin{equation}
    f_{21}(t):=-(\mathbbm1-\mathcal C_r(t))^{-1}\mathcal C(t)P_b(t)P^+e^{-itH_0}f,
\end{equation}
\begin{equation}
    f_{22}(t):=-(\mathbbm1-\mathcal C_r(t))^{-1}\mathcal C_M(t)P^+e^{-itH_0}f
\end{equation}
and
\begin{equation}
    f_{23}(t):=(\mathbbm1-\mathcal C_r(t))^{-1}\mathcal C(t)P^+e^{-itH_0}f.
\end{equation}
By Proposition~\ref{prop: weight}, Corollary~\ref{cor: bd: C}, estimate~\eqref{Nov.20.6} and the $L^2$ local decay of the free flow, $f_{21}(t)$ satisfies for all $n\geq 8$ and $\eta\geq 0$,
\begin{align}
  \left( \int_0^\infty \| \langle x\rangle^{-\eta}f_{21}(t)\|^2dt\right)^{1/2}\lesssim &\left( \int_0^\infty \|P_b(t) P^+e^{-itH_0}f\|^2 dt\right)^{1/2} \nonumber\\
  \lesssim & \left( \int_0^\infty \|P_b(t)\langle x\rangle^{3/2}\|\|\langle x\rangle^{-3/2} P^+ \langle x\rangle^{3/2}\|\|\langle x\rangle^{-3/2} e^{-itH_0}f\|^2 dt\right)^{1/2}\nonumber\\
  \lesssim & \|f\|.
\end{align}
By Proposition~\ref{prop: 4.1}, estimates~\eqref{BP2: est1} and~\eqref{Aug.8} and $L^2$ local smoothing estimate of the free flow, $f_{22}(t)$ satisfies 
\begin{align}
    & \left( \int_0^\infty \| \langle x\rangle^{-\eta}f_{22}(t)\|^2dt\right)^{1/2}\nonumber\\
    \lesssim &\left( \int_0^\infty \|\mathcal C_M(t) P^+e^{-itH_0}f\|^2 dt\right)^{1/2} \nonumber\\
  \lesssim & \left( \int_0^\infty \|\mathcal C_M(t)|p|^{-1/2}\langle x\rangle^{1/2+\epsilon}\|\|\langle x\rangle^{-1/2-\epsilon} |p|^{1/2}P^+ |p|^{-1/2}\langle x\rangle^{1/2+\epsilon}\|\|\langle x\rangle^{-1/2-\epsilon} |p|^{1/2} e^{-itH_0}f\|^2 dt\right)^{1/2}\nonumber\\
  \lesssim_{M,\epsilon} & \|f\|
\end{align}
for all $\epsilon\in (0,\frac{1}{2})$. Next, we estimate~$\langle x\rangle^{-\eta}f_{23}(t)$. For this, we note that by Proposition~\ref{prop: main2} and Eqs.~\eqref{def: mCt} and~\eqref{def: Ipmt}, for all $\alpha\in [0,\frac{7}{4})$ and $n\geq8$, with $(\alpha-\frac{1}{8})_+:=\max\{0,\alpha-\frac{1}{8} \}$,
\begin{align}
\| \frac{|p|^{\alpha-\frac{1}{8}}}{\langle p\rangle^{\alpha-\frac{1}{8}}}\mathcal C(t)\frac{\langle p\rangle^\alpha}{|p|^\alpha}\|\lesssim_n 1.
\end{align}
This together with Corollary~\ref{cor: bd: C}, Proposition~\ref{prop: weight} and estimates~\eqref{BP2: est2} and $\|\frac{|p|^\alpha}{\langle p\rangle^\alpha}\|\leq 1$ implies 
\begin{align}
\| \frac{|p|^{(\alpha-\frac{1}{8})_+}}{\langle p\rangle^{(\alpha-\frac{1}{8})_+}}\mathcal C_r(t)\frac{\langle p\rangle^\alpha}{|p|^\alpha}\|\leq & \|\frac{|p|^{(\alpha-\frac{1}{8})_+}}{\langle p\rangle^{(\alpha-\frac{1}{8})_+}}\mathcal C(t)\frac{\langle p\rangle^\alpha}{|p|^\alpha}\|+\|\frac{|p|^{(\alpha-\frac{1}{8})_+}}{\langle p\rangle^{(\alpha-\frac{1}{8})_+}}\mathcal C_M(t)\frac{\langle p\rangle^\alpha}{|p|^\alpha}\|+\|\frac{|p|^{(\alpha-\frac{1}{8})_+}}{\langle p\rangle^{(\alpha-\frac{1}{8})_+}}\mathcal C(t)P_b(t)\frac{\langle p\rangle^\alpha}{|p|^\alpha}\|\nonumber\\
\lesssim_{n}&1+\|C_M(t)\frac{\langle p\rangle^\alpha}{|p|^\alpha} \|+\| P_b(t)\langle x\rangle^\alpha\|\| \langle x\rangle^{-\alpha}\frac{\langle p\rangle^\alpha}{|p|^\alpha}\|\nonumber\\
\lesssim_{n,M,\alpha}& 1,\label{est: iter}
\end{align}
where we also used Hardy-Littlewood-Sobolev inequality. This together with Hardy-Littlewood-Sobolev inequality and Neumann series yields for all $\eta>\frac{3}{2}$ and $\alpha\in (\frac{3}{2}, \min\{\eta, \frac{7}{4}\})$, 
\begin{equation}
    \| \langle x\rangle^{-\eta}(\mathbbm 1-\mathcal C_r(t))^{-1}\frac{\langle p\rangle^
\alpha}{|p|^\alpha}\|\lesssim_{n,M,\eta,\alpha} 1. 
\end{equation}
We also note that Proposition~\ref{prop: main} and  Eq.~\eqref{def: mCt} yields 
\begin{equation}
   \sup\limits_{t\in \mathbb R} \|\frac{|p|^\alpha}{\langle p\rangle^\alpha} \mathcal C(t)|p|^{-1/2}\langle x \rangle^{1/2+\epsilon}\|\lesssim_{\epsilon,n}1
\end{equation}
holds true for all $n\geq 5$, $\epsilon\in (0, \frac{1}{2})$ and $\alpha\in (\frac{3}{2}, \frac{7}{4})$. This together with Corollary~\ref{cor: bd: C} and estimates $\|\frac{|p|^\alpha}{\langle p\rangle^\alpha}\|\leq 1, \alpha\geq 0$ and~\eqref{BP2: est1}, yields 
\begin{align}
    \sup\limits_{t\in \mathbb R} \|\frac{|p|^\alpha}{\langle p\rangle^\alpha} \mathcal C_r(t)|p|^{-1/2}\langle x \rangle^{1/2+\epsilon}\|\leq & \sup\limits_{t\in \mathbb R} \|\frac{|p|^\alpha}{\langle p\rangle^\alpha} \mathcal C(t)|p|^{-1/2}\langle x \rangle^{1/2+\epsilon}\|+\sup\limits_{t\in \mathbb R} \|\frac{|p|^\alpha}{\langle p\rangle^\alpha} \mathcal C_M(t)|p|^{-1/2}\langle x \rangle^{1/2+\epsilon}\|\nonumber\\
    &+\sup\limits_{t\in \mathbb R} \|\frac{|p|^\alpha}{\langle p\rangle^\alpha} \mathcal C(t)P_b(t)|p|^{-1/2}\langle x \rangle^{1/2+\epsilon}\|\nonumber\\
    \lesssim_{\epsilon,n} & 1+  \sup\limits_{t\in \mathbb R} \| \mathcal C_M(t)|p|^{-1/2}\langle x \rangle^{1/2+\epsilon}\|+\sup\limits_{t\in \mathbb R} \|P_b(t)|p|^{-1/2}\langle x \rangle^{1/2+\epsilon}\|\nonumber\\
    \lesssim_{\epsilon, n, M}&1+\sup\limits_{t\in \mathbb R} \|P_b(t)|p|^{-1/2}\langle x \rangle^{1/2+\epsilon}\|.
\end{align}
This together with Proposition~\ref{prop: weight} and estimate~\eqref{Blem: est} yields for $n\geq 8$,
\begin{align}
     \sup\limits_{t\in \mathbb R} \|\frac{|p|^\alpha}{\langle p\rangle^\alpha} \mathcal C_r(t)|p|^{-1/2}\langle x \rangle^{1/2+\epsilon}\|\lesssim_{\epsilon, n, M}&1+\sup\limits_{t\in \mathbb R} \|P_b(t)\langle x\rangle^{\frac{3}{2}}\|\|\langle x\rangle^{-3/2} |p|^{-1/2}\langle x \rangle^{1/2+\epsilon}\|\nonumber\\
     \lesssim_{\epsilon, n, M}& 1.
\end{align}
We fix $M=M_0$. This together with estimates~\eqref{est: iter} and~\eqref{Aug.8} and Neumann series yields 
\begin{equation}
    \sup\limits_{t\in \mathbb R} \|\langle x\rangle^{-\eta}(\mathbbm 1-\mathcal C_r(t))^{-1}P^+|p|^{-1/2}\langle x\rangle^{1/2+\epsilon}\|\lesssim_{n,\eta}1
\end{equation}
for all $\epsilon\in (0, \frac{1}{2})$ and $\eta> \frac{3}{2}$. Therefore, by the Restricted $L^2$ local decay of the free flow, we conclude 
\begin{equation}
    \left(\int_0^\infty  \| \langle x\rangle^{-\eta} (\mathbbm1 -\mathcal C_r(t))^{-1} P^+ e^{-itH_0}f\|^2dt\right)^{1/2}\lesssim_{n,\eta} \|f\|.
\end{equation}
Similarly, we have 
\begin{equation}
    \left(\int_0^\infty  \| \langle x\rangle^{-\eta} (\mathbbm1 -\mathcal C_r(t))^{-1} P^- e^{-itH_0}f\|^2dt\right)^{1/2}\lesssim_{n,\eta} \|f\|.
\end{equation}\end{proof}
Now we prove Proposition~\ref{BP2}. 
\begin{proof}[Proof of Proposition~\ref{BP2}] Using Assumption~\ref{aspV3} instead of Proposition~\ref{prop: weight} and following a similar proof of Proposition~\ref{BP1}, we obtain~\eqref{est: BP2}.
    
\end{proof}

\subsection{Proof of Theorems~\ref{thm} and~\ref{thm2}}
Now we prove Theorems~\ref{thm} and~\ref{thm2}.
\begin{proof}[Proof of Theorem~\ref{thm}] By Proposition~\ref{prop: 4.1}, we obtain that there exists $M_0>0$ such that whenever $M\geq M_0$,
\begin{align}
    U(t,0)f=&(\mathbbm1-\mathcal C_r(t))^{-1}P^+e^{-itH_0}\Omega_+^*(0)f+(\mathbbm1-\mathcal C_r(t))^{-1}P^-e^{-itH_0}\Omega_-^*(0)f\nonumber\\
    &+(\mathbbm1-\mathcal C_r(t))^{-1}\mathcal C_{M}(t)U(t,0)f
\end{align}
for all $f\in \mathcal L^2_x(\mathbb R^n)$ with $f=P_c(0)f$ and $\| (\mathbbm1-\mathcal C_r(t))^{-1}\|\leq 2$ for all $t\in \mathbb R$. This together with Propositions~\ref{BP1} and estimates~\eqref{id: CMt} and 
\begin{equation}
    \left(\int_0^\infty \|F_M(x,p)e^{-itH_0}f\|^2dt\right)^{1/2}\lesssim_M\|f\|,\qquad \forall\, f\in \mathcal L^2_x(\mathbb R^n)
\end{equation}
yields
\begin{align}
  \left(\int_0^\infty  \| \langle x\rangle^{-\eta}U(t,0)f\|^2dt\right)^{1/2}\leq & \left(\int_0^\infty  \| \langle x\rangle^{-\eta}(\mathbbm1-\mathcal C_r(t))^{-1}P^+e^{-itH_0}\Omega_+^*(0)f\|^2dt\right)^{1/2}\nonumber\\
  &+\left(\int_0^\infty  \| \langle x\rangle^{-\eta}(\mathbbm1-\mathcal C_r(t))^{-1}P^-e^{-itH_0}\Omega_+^*(0)f\|^2dt\right)^{1/2}\nonumber\\
  & +\left(\int_0^\infty  \| \langle x\rangle^{-\eta}(\mathbbm1-\mathcal C_r(t))^{-1}\mathcal C_M(t)U(t,0)f\|^2dt\right)^{1/2}\nonumber\\
  \lesssim_{\eta,n, M} &\|f\|
\end{align}
for all $\eta>\frac{3}{2}$ and $f\in \mathcal L^2_x(\mathbb R^n)$ with $f=P_c(0)f$ and $n\geq 8$.
\end{proof}
\begin{proof}[Proof of Theorem~\ref{thm2}] Following a similar argument of Theorem~\ref{thm} and using Proposition~\ref{BP2} instead of Proposition~\ref{BP1}, we arrive at 
 \begin{align}
  \left(\int_0^\infty  \| \langle x\rangle^{-\eta}U(t,0)f\|^2dt\right)^{1/2}\leq & \left(\int_0^\infty  \| \langle x\rangle^{-\eta}(\mathbbm1-\mathcal C_r(t))^{-1}P^+e^{-itH_0}\Omega_+^*(0)f\|^2dt\right)^{1/2}\nonumber\\
  &+\left(\int_0^\infty  \| \langle x\rangle^{-\eta}(\mathbbm1-\mathcal C_r(t))^{-1}P^-e^{-itH_0}\Omega_+^*(0)f\|^2dt\right)^{1/2}\nonumber\\
  & +\left(\int_0^\infty  \| \langle x\rangle^{-\eta}(\mathbbm1-\mathcal C_r(t))^{-1}\mathcal C_M(t)U(t,0)f\|^2dt\right)^{1/2}\nonumber\\
  \lesssim_{\eta,n, M} &\|f\|
\end{align}
for all $\eta>\frac{3}{2}$ and $f\in \mathcal L^2_x(\mathbb R^n)$ with $f=P_c(0)f$ and $n\geq 5$.
\end{proof}

\subsection{Strichartz estimates}
The Strichartz estimates for the free flow (see \cite{S1977}, \cite{T2006} for Strichartz estimates and \cite{KT1998} for endpoint Strichartz estimates) state that for \(2 \leq r, q \leq \infty\), \(\frac{n}{r} + \frac{2}{q} = \frac{n}{2}\), and \((q, r, n) \neq (2, \infty, 2)\), the following estimates hold:

\begin{enumerate}
    \item \textbf{Homogeneous Strichartz Estimate}:
    \[
    \| e^{-itH_0} f \|_{\s^q_t \s^r_x(\mathbb{R}^{n+1})} \lesssim_{n, q, r} \| f \|_{\s^2_x(\mathbb{R}^n)}
    \]

    \item \textbf{Dual Homogeneous Strichartz Estimate}:
    \[
    \left\| \int_\mathbb{R} dse^{-isH_0}F(s) \right\|_{\s^2_x(\mathbb{R}^n)} \lesssim_{n, q, r} \| F \|_{\s^{q'}_t \s^2_x(\mathbb{R}^n)}
    \]

    \item \textbf{Inhomogeneous Strichartz Estimate}:
    \[
    \left\| \int_{s<t} dse^{-i(t-s)H_0}F(s) \right\|_{\s^q_t \s^r_x(\mathbb{R}^{n+1})} \lesssim_{n, q, r} \| F \|_{\s^{q'}_t \s^{r'}_x(\mathbb{R}^{n+1})}
    \]
\end{enumerate}
where \((r, r')\) and \((q, q')\) are conjugate pairs.

For a perturbed system, the corresponding Strichartz estimates are
\[
\| U(t,0)P_c(0) f \|_{\s^q_t \s^r_x(\mathbb{R}^{n+1})} \leq C_q \| f \|_{\s^2_x(\mathbb{R}^n)}
\]
for \(2 \leq r, q \leq \infty\), \(\frac{n}{r} + \frac{2}{q} = \frac{n}{2}\), and \((q, r, n) \neq (2, \infty, 2)\). In this subsection, we prove that our local decay estimates imply Strichartz estimates by utilizing the inhomogeneous Strichartz estimate for the free flow.

\begin{proof}[Proof of Theorem \ref{thmS}]It suffices to check the endpoint Strichartz estimates, that is, the case when $(q,r,n)=(2,\frac{2n}{n-2},n)$, $n\geq 5$. By the Duhamel's principle,
\begin{align}
U(t,0)P_c(0)\psi=& e^{-itH_0}\psi+(-i)\int_0^t ds e^{-i(t-s)H_0}V(x,s)U(s,0)P_c(0)\psi\\
=:&\psi_1(t)+\psi_2(t).
\end{align}
$\psi_1(t)$ enjoys Strichartz estimates  of the free flow. For $\psi_2(t)$, by the inhomogeneous Strichartz estimate of the free flow, the bound 
\begin{equation}
    \|\langle x\rangle^{-2}\|_{\mathcal L^2_x(\mathbb R^n)\to \mathcal L^{\frac{2n}{n+2}}_x(\mathbb R^{n})}\lesssim 1
\end{equation}
and Theorem~\ref{thm2}, we have
\begin{align}
\|\psi_2(t) \|_{\s^2_t\s^{\frac{2n}{n+2}}_x(\mathbb{R}^{n}\times \mathbb R^+)}\lesssim &\|  V(x,t)U(t,0)P_c(0)\psi \|_{\s^{2}_t\s^{\frac{2n}{n+2}}_x(\mathbb{R}^{n}\times \mathbb R^+)}\\
\lesssim & \| \langle x\rangle^{2} V(x,t)\|_{\s^\infty_{t}\s^\infty_x(\mathbb{R}^{n+1})}\| \langle x\rangle^{-2}U(t,0)P_c(0)\psi \|_{\s^2_{x,t}(\mathbb{R}^{n}\times \mathbb R^+)}\\
\lesssim & \| \langle x\rangle^{2}V(x,t)\|_{\s^\infty_{x,t}(\mathbb{R}^{n+1})}\| \psi\|_{\s^2_x(\mathbb{R}^n)}.
\end{align}
This together with the estimate for $\psi_1(t)$ yields 
\eq
\|\psi(t) \|_{\s^2_t\s^{\frac{2n}{n+2}}_x(\mathbb{R}^{n}\times \mathbb R^+)}\lesssim\| \langle x\rangle^{2}V(x,t)\|_{\s^\infty_{x,t}(\mathbb{R}^{n+1})}\| \psi\|_{\s^2_x(\mathbb{R}^n)}.
\eeq\end{proof}
\appendix
\section{Proof of free estimates }\label{sec: app}
\begin{proof}[Proof of Lemma~\ref{out/in1}]\textbf{Proof of \eqref{Sep29.1}:} It suffices to check the case when $t>1$. Let
\eq
F_{2^j}(|p|):= F(2^j\leq |p|<   2^{j+1}):=F(|p|\geq 2^j)-F(|p|<2^{j+1}),\quad j=1,2,\cdots.
\eeq
Then 
\eq
\sum\limits_{j=1}^\infty F_{2^j}(|p|)=F(|p|\geq 2).
\eeq
By estimate~\eqref{est: Mourre}, it suffices to prove
\eq
\| \sum\limits_{j=1}^\infty P^{\pm}  F_{2^j}(|p|)e^{\pm i tH_0}\langle x\rangle^{-\delta}\|\lesssim_{n}  \frac{1}{\langle t\rangle^{\delta}}.\label{A5}
\eeq
Indeed, using dilation transformation, we have
\begin{multline}
\| P^{\pm}  F_{2^j}(|p|)e^{\pm i tH_0}\langle x\rangle^{-\delta}\|=\| P^{\pm}  F_{2}(|p|)e^{\pm i 2^{2j}tH_0}\langle x/2^j\rangle^{-\delta}\|\\
\leq \| P^{\pm}  F_{2}(|p|)e^{\pm i 2^{2j}tH_0}\langle x\rangle^{-\delta}\|\times\| \langle x\rangle^\delta\langle x/2^j\rangle^{-\delta}\|\\
\lesssim_n \frac{1}{\langle 2^{2j}t\rangle^\delta} \times 2^{j\delta}\lesssim_n \frac{1}{2^{j\delta}} \times \frac{1}{\langle t\rangle^{\delta}},\label{A4}
\end{multline}
which implies~\eqref{A5}.\\
\textbf{Proof of \eqref{Sep20.2}: }Using the dilation transformation, 
\eq
\| P^{\pm}  F(|p|\geq M)e^{\pm i tH_0}|p|^l\langle x\rangle^{-\delta}\|=M^l\| P^{\pm}  F(|p|>1)e^{\pm i M^2tH_0}|p|^l\langle x/M\rangle^{-\delta}\|.
\eeq
Proceeding as~\eqref{A4}, we have for $M^2t\geq 1$, 
\begin{align}
    \| P^{\pm}  F(|p|\geq M)e^{\pm i tH_0}|p|^l\langle x\rangle^{-\delta}\|\leq &\sum\limits_{j=0}^\infty M^l\| P^{\pm}  F_{2^j}(|p|)e^{\pm i M^2tH_0}|p|^l\langle x/M\rangle^{-\delta}\|\\
    \lesssim_n &\sum\limits_{j=0}^\infty  2^{jl}M^l\times \frac{1}{2^{j\delta }\langle M^2t\rangle^\delta}M^\delta\\
    \lesssim_{n, l-\delta} &\frac{1}{M^{\delta-l}t^\delta}. 
\end{align}
\textbf{Proof of \eqref{Oct.1}: } We note that
\begin{multline}
\int_0^1t^2dt\| P^{\pm}  F_{2^j}(|p|)e^{\pm i tH_0}|p|^2\langle x\rangle^{-\delta}\|\leq \int_0^{1/2^{3j/4}}t^2dt\| P^{\pm}  F_{2^j}(|p|)e^{\pm i tH_0}|p|^2\langle x\rangle^{-\delta}\|+\\
\int_{1/2^{3j/4}}^1t^2dt\| P^{\pm}  F_{2^j}(|p|)e^{\pm i tH_0}|p|^2\langle x\rangle^{-\delta}\|=:A_{j,1}+A_{j,2}.\label{def: Aj1Aj2}
\end{multline}
For $A_{j,1}$, we have
\eq
A_{j,1}\lesssim \int_0^{1/2^{3j/4}}dt t^2 \times 2^{2j}\lesssim 2^{-j/4}.\label{est: Aj1}
\eeq
For $A_{j,2}$, using dilation to replace $|p|$ with $2^j|p|$, we have
\begin{align}
A_{j,2}=&2^{2j}\int_{1/2^{3j/4}}^1 t^2dt \| P^{\pm}  F_{1}(|p|)e^{\pm i 2^{2j}tH_0}|p|^2\langle x/2^j\rangle^{-\delta}\|\\
\lesssim_{n,\delta}& 2^{2j}\int_{1/2^{3j/4}}^1 t^2dt\times \frac{1}{\langle t2^{2j}\rangle^\delta}\times 2^{j\delta}\\
\lesssim_{n,\delta}& \frac{1}{2^{j(\delta-2)}}.
\end{align}
This together with estimate~\eqref{est: Aj1} and Eq.~\eqref{def: Aj1Aj2} yields 
\begin{equation}
    \begin{aligned}
        \int_0^1t^2dt\| P^{\pm}  F(|p|\geq 1)e^{\pm i tH_0}|p|^2\langle x\rangle^{-\delta}\|\leq &\sum\limits_{j=0}^\infty \int_0^1t^2dt\| P^{\pm}  F_{2^j}(|p|)e^{\pm i tH_0}|p|^2\langle x\rangle^{-\delta}\|\\
        \lesssim_{n,\delta} & \sum\limits_{j=0}^\infty \left( 2^{-j/4}+2^{-j(\delta-2)}\right)\\
         \lesssim_{n,\delta}& 1.
    \end{aligned}
\end{equation}
\textbf{Proof of \eqref{Sep20.1}: } It suffices to estimate~\eqref{Sep20.1} when $t\geq 1$. By estimate~\eqref{Sep20.2} with $l=0$ and $M=\frac{1}{\langle t\rangle^{1/2-\epsilon}}$, we have 
\begin{equation}
\| \frac{|p|^\alpha}{\langle p\rangle^\alpha }P^{\pm}  F(|p|\geq\frac{1}{\langle t\rangle^{1/2-\epsilon}})e^{\pm i tH_0} \langle x\rangle^{-\delta}\|\lesssim_{n, \delta} \frac{1}{\langle t\rangle^{(\frac{1}{2}+\epsilon)\delta}}.\label{est: highA}
\end{equation}
We also note that by estimate~\eqref{paf}, we have 
\begin{equation}
    \begin{aligned}
        &\| \frac{|p|^\alpha}{\langle p\rangle^\alpha }P^{\pm}  F(|p|< \frac{1}{\langle t\rangle^{1/2-\epsilon}})e^{\pm i tH_0} \langle x\rangle^{-\delta}\|\\
        \lesssim_{\alpha, R} & \| |p|^\alpha F(|p|< \frac{1}{\langle t\rangle^{1/2-\epsilon}})e^{\pm i tH_0} \langle x\rangle^{-\delta}\| \\
        \lesssim_{\alpha,R} & \| |p|^{\alpha+\min\{\frac{n}{2}-\epsilon,\delta\}} F(|p|< \frac{1}{\langle t\rangle^{1/2-\epsilon}})\|\| e^{\pm i tH_0} \|\||p|^{-\min\{\frac{n}{2}-\epsilon,\delta\}}\langle x\rangle^{-\delta}\|\\
        \lesssim_{\alpha,R,n} & \frac{1}{\langle t\rangle^{(1/2-\epsilon)(\alpha+\min\{\frac{n}{2}-\epsilon,\delta\})}}.
    \end{aligned}\label{est: lowA}
\end{equation}
Estimates~\eqref{est: highA} and~\eqref{est: lowA} yield~\eqref{Sep20.1}.\\
\textbf{Proof of~\eqref{est: Pbweight}:} We note that similarly, for $\epsilon\in (0,1/2)$,  
\begin{equation}
    \|  \frac{1}{\langle x\rangle^{2-\delta}}P^\pm e^{\pm itH_0}F(|p|\geq\frac{1}{\langle t\rangle^{1/2-\epsilon}})p_j^l \langle x\rangle^{-(\frac{n}{2}+l+2-\delta)}\|\lesssim_{n,\delta} \frac{1}{\langle t\rangle^{\frac{n}{4}+\frac{l+2-\delta}{2}+\epsilon}},\qquad l=1,2
\end{equation}
and 
\begin{equation}
\begin{aligned}
   & \|  \frac{1}{\langle x\rangle^{2-\delta}}P^\pm e^{\pm itH_0}F(|p|< \frac{1}{\langle t\rangle^{1/2-\epsilon}})p_j^l \langle x\rangle^{-(\frac{n}{2}+l+2-\delta)}\|\\
   \lesssim_{n,\delta, R} & \|  e^{\pm itH_0}F(|p|< \frac{1}{\langle t\rangle^{1/2-\epsilon}})|p|^{2-\delta+n/2-\epsilon}p_j^l \|\\
   \lesssim_{n,\delta,R}& \frac{1}{\langle t\rangle^{(n/2+2+l-\delta-\epsilon)(1/2-\epsilon)}}.
\end{aligned}
\end{equation}
These yield~\eqref{est: Pbweight}.

\end{proof}
\begin{proof}[Proof of Lemma~\ref{lem: cpt F}] Let $\bar F_{c,1}(z) := 1 - F_{c,1}(z)$ denote the complement of $F_{c,1}(z)$. Then 
\begin{equation}
    F_{c,1}(H)-F_{c,1}(H_0)=\bar F_{c,1}(H_0)-\bar F_{c,1}(H). 
\end{equation}
Using the Fourier inversion Theorem, we have 
\begin{equation}
    \bar F_{c,1}(H_0)-\bar F_{c,1}(H)=\frac{1}{\sqrt{2\pi}}\int \hat F_{c,1}(w) \left( e^{iwH_0}-e^{iwH}\right)dw
\end{equation}
where $\hat F_{c,1}(w)$ denotes the Fourier transform of $F_{c,1}(z)$. By the Duhamel's principle, this implies 
\begin{equation}
     \bar F_{c,1}(H_0)-\bar F_{c,1}(H)=\frac{-i}{\sqrt{2\pi}}\int \hat F_{c,1}(w) \int_0^w e^{i(w-u)H_0}Ve^{iuH}du dw.
\end{equation}
Since $ \langle w\rangle\hat F_{c,1}(w)\in \mathcal L^1_w(\mathbb R)$, it suffices to show that 
\begin{equation}
    \int_0^w e^{i(w-u)H_0}Ve^{iuH} du\text{ is compact for all $w\in \mathbb R$.}
\end{equation}
Since $F(|p|\leq 1)\langle x\rangle^{-1} $ is compact, it suffices to show that 
\begin{equation}
   \mathcal O:= \int_0^w F(|p|\geq 1)e^{i(w-u)H_0}Ve^{iuH} du\text{ is compact for all $w\in \mathbb R$.}
\end{equation}
For this, we write $\mathcal O$ as 
\begin{equation}
    \mathcal O=\int_0^w |p|^{1/2}e^{i(w-u)H_0}\langle x\rangle^{-1}\langle x\rangle|p|^{-1/2}F(|p|\geq1)Ve^{iuH} du.
\end{equation}
We note that 
\begin{equation}
    \langle x\rangle|p|^{-1/2}F(|p|\geq1)V=|p|^{-1/2}F(|p|\geq1)\langle x\rangle V+[\langle x\rangle, |p|^{-1/2}F(|p|\geq1)]V
\end{equation}
and therefore, $ \langle x\rangle|p|^{-1/2}F(|p|\geq1)V$ is compact. This together with the Restricted $L^2$ local decay estimate yields the compactness of $\mathcal O$. \end{proof}
\begin{proof}[Proof of Lemma~\ref{lem: Fwell}] We estimate $F_{c,1}(H)(\mathbbm1-F_{c,2}(H_0))|p|^{-1/2}\langle x\rangle^{1/2+\epsilon}$ first. Let $\bar F_{c,j}:=1-F_{c,j}, j=1,2.$ Using the Fourier inversion Theorem and the Duhamel's principle, we obtain 
\begin{align}
    F_{c,1}(H)\bar F_{c,2}(H_0)=(-i)\int \int_0^w\hat{\bar{F}}_{c,2}(w) F_{c,1}(H)e^{i(w-s)H} Ve^{isH_0}dsdw,
\end{align}
where we also used $F_{c,1}(H)\bar F_{c,2}(H)=0$. This together with~\eqref{B21} yields 
\begin{align}
    \| F_{c,1}(H)\bar F_{c,2}(H_0)|p|^{-1/2}\langle x\rangle^{1/2+\epsilon}\|\lesssim & \int \langle w\rangle^{1+2\epsilon}|\hat{\bar{F}}_{c,2}(w)|\|\langle x\rangle^2V(x)\|_{\mathcal L^\infty_x(\mathbb R^3)} \nonumber\\
    &\times\int  \langle s\rangle^{-1-2\epsilon}\| \langle x\rangle^{-2}e^{-isH_0}|p|^{-1/2}\langle x\rangle^{1/2+\epsilon}\|dsdw\nonumber\\
    \lesssim_\epsilon &  \|\langle x\rangle^2V(x)\|_{\mathcal L^\infty_x(\mathbb R^3)} .
\end{align}
Similarly, we have~\eqref{lem3.2: est: 2}.\end{proof}

\begin{proof}[Proof of Lemma~\ref{lem: cpt F}] Let $\bar F_{c,1}(z) := 1 - F_{c,1}(z)$ denote the complement of $F_{c,1}(z)$. Then 
\begin{equation}
    F_{c,1}(H)-F_{c,1}(H_0)=\bar F_{c,1}(H_0)-\bar F_{c,1}(H). 
\end{equation}
Using the Fourier inversion Theorem, we have 
\begin{equation}
    \bar F_{c,1}(H_0)-\bar F_{c,1}(H)=\frac{1}{\sqrt{2\pi}}\int \hat F_{c,1}(w) \left( e^{iwH_0}-e^{iwH}\right)dw
\end{equation}
where $\hat F_{c,1}(w)$ denotes the Fourier transform of $F_{c,1}(z)$. By the Duhamel's principle, this implies 
\begin{equation}
     \bar F_{c,1}(H_0)-\bar F_{c,1}(H)=\frac{-i}{\sqrt{2\pi}}\int \hat F_{c,1}(w) \int_0^w e^{i(w-u)H_0}Ve^{iuH}du dw.
\end{equation}
Since $ \langle w\rangle\hat F_{c,1}(w)\in \mathcal L^1_w(\mathbb R)$, it suffices to show that 
\begin{equation}
    \int_0^w e^{i(w-u)H_0}Ve^{iuH} du\text{ is compact for all $w\in \mathbb R$.}
\end{equation}
Since $F(|p|\leq 1)\langle x\rangle^{-1} $ is compact, it suffices to show that 
\begin{equation}
   \mathcal O:= \int_0^w F(|p|\geq 1)e^{i(w-u)H_0}Ve^{iuH} du\text{ is compact for all $w\in \mathbb R$.}
\end{equation}
For this, we write $\mathcal O$ as 
\begin{equation}
    \mathcal O=\int_0^w |p|^{1/2}e^{i(w-u)H_0}\langle x\rangle^{-1}\langle x\rangle|p|^{-1/2}F(|p|\geq1)Ve^{iuH} du.
\end{equation}
We note that 
\begin{equation}
    \langle x\rangle|p|^{-1/2}F(|p|\geq1)V=|p|^{-1/2}F(|p|\geq1)\langle x\rangle V+[\langle x\rangle, |p|^{-1/2}F(|p|\geq1)]V
\end{equation}
and therefore, $ \langle x\rangle|p|^{-1/2}F(|p|\geq1)V$ is compact. This together with the Restricted $L^2$ local decay estimate yields the compactness of $\mathcal O$. \end{proof}
\begin{proof}[Proof of Lemma~\ref{lem: Fwell}] We estimate $F_{c,1}(H)(\mathbbm1-F_{c,2}(H_0))|p|^{-1/2}\langle x\rangle^{1/2+\epsilon}$ first. Let $\bar F_{c,j}:=1-F_{c,j}, j=1,2.$ Using the Fourier inversion Theorem and the Duhamel's principle, we obtain 
\begin{align}
    F_{c,1}(H)\bar F_{c,2}(H_0)=(-i)\int \int_0^w\hat{\bar{F}}_{c,2}(w) F_{c,1}(H)e^{i(w-s)H} Ve^{isH_0}dsdw,
\end{align}
where we also used $F_{c,1}(H)\bar F_{c,2}(H)=0$. This together with~\eqref{B21} yields 
\begin{align}
    \| F_{c,1}(H)\bar F_{c,2}(H_0)|p|^{-1/2}\langle x\rangle^{1/2+\epsilon}\|\lesssim & \int \langle w\rangle^{1+2\epsilon}|\hat{\bar{F}}_{c,2}(w)|\|\langle x\rangle^2V(x)\|_{\mathcal L^\infty_x(\mathbb R^3)} \nonumber\\
    &\times\int  \langle s\rangle^{-1-2\epsilon}\| \langle x\rangle^{-2}e^{-isH_0}|p|^{-1/2}\langle x\rangle^{1/2+\epsilon}\|dsdw\nonumber\\
    \lesssim_\epsilon &  \|\langle x\rangle^2V(x)\|_{\mathcal L^\infty_x(\mathbb R^3)} .
\end{align}
Similarly, we have~\eqref{lem3.2: est: 2}.\end{proof}

\bigskip
\paragraph{\bf Acknowledgment}

A.S. is partially supported by NSF-DMS-220931. X.W. is partially supported by ARC-FL220100072, NSF-DMS-220931 and NSERC Grant NA7901.  The authors thank Maxime Van de Moortel for careful reading and for useful discussions.

Parts of this work were done while the third author was at the Fields Institute  for Research in Mathematical Sciences, Toronto, Texas A\&M University, Rutgers University and University of Toronto. 

\bibliographystyle{abbrv}

\end{document}